\documentclass[11pt,a4paper]{amsart}

\usepackage{amssymb}
\usepackage{graphicx}
\usepackage{pst-all}
\usepackage{url}

\newtheorem{theorem}{Theorem}[section]

\theoremstyle{definition}

\theoremstyle{remark}
\newtheorem{remark}[theorem]{Remark}

\numberwithin{equation}{section}

\newcommand{\softd}{{\leavevmode\setbox1=\hbox{d}%
          \hbox to 1.05\wd1{d\kern-0.4ex{\char039}\hss}}}
\newcommand{\abs}[1]{\lvert#1\rvert}
\newcommand{\norm}[1]{\lVert#1\rVert}
\newcommand{\D}{\partial}

\newcommand{\Dt}{\partial t}
\newcommand{\dt}{\Delta t}
\newcommand{\dd}{\mathrm{d}}

\def\bold#1{\mbox{\boldmath $#1$}}
\newcommand{\uu}[1]{\bold{#1}}

\newcommand{\const}{\mathrm{const}}

\begin{document}

\title[Propagation of a 3-D Weak Shock Front]{Propagation of a
  Three-dimensional Weak Shock Front Using Kinematical Conservation
  Laws}

\author[Arun]{K. R. Arun}
 \thanks{K. R. A. gratefully acknowledges the financial support from the
   Alexander-von-Humboldt Foundation through a postdoctoral fellowship
   during early stages of this work.}
\address{School of Mathematics, Indian Institute of Science Education and 
  Research, Thiruvananthapuram - 695016, India.}
\email{arun@iisertvm.ac.in}

\author[Prasad]{Phoolan Prasad}
\thanks{This research work was partly completed when P. P. was a Raja
  Ramanna Fellow, supported by the Department of Atomic Energy,
  Government of India. The work was completed when this author was
  supported by the National Academy of Sciences, India, through the
  NASI-Senior Scientist Platinum Jubilee Fellowship. The Department of
  Mathematics, Indian Institute of Science, was supported the
  University Grants Commission, India, through the UGC-SAP-Centre for
  Advanced Study and the Department of Science and Technology,
  Government of India, through the FIST Programme.}
\address{Department of Mathematics, Indian Institute of Science,
  Bangalore - 560012, India.}
\email{prasad@math.iisc.ernet.in}
\urladdr{http://math.iisc.ernet.in/~prasad/}


\date{\today}


\begin{abstract}
  In this paper we present a mathematical theory and a numerical
  method to study the propagation of a three-dimensional (3-D) weak
  shock front into a polytropic gas in a uniform state and at rest,
  though the method can be extended to shocks moving into nonuniform
  flows. The theory is based on the use of 3-D kinematical
  conservation laws (KCL), which govern the evolution of a surface in
  general and a shock front in particular. The 3-D KCL, derived purely
  on geometrical considerations, form an under-determined system of
  conservation laws. In the present paper the 3-D KCL system is closed
  by using two appropriately truncated transport equations from an
  infinite hierarchy of compatibility conditions along shock rays. The
  resulting governing equations of this KCL based 3-D shock ray
  theory, leads to a weakly hyperbolic system of eight conservation
  laws with three divergence-free constraints. The conservation laws
  are solved using a Godunov-type central finite volume scheme, with a
  constrained transport technique to enforce the constraints. The
  results of extensive numerical simulations reveal several physically
  realistic geometrical features of shock fronts and the complex
  structures of kink lines formed on them. A comparison of the results
  with those of a weakly nonlinear wavefront shows that a weak shock
  front and a weakly nonlinear wavefront are topologically same. The
  major important differences between the two are highlighted in the
  contexts of corrugational stability and converging shock fronts.
\end{abstract}

\newpage

\maketitle

\section{Introduction}
\label{sect:intro}

A very attractive method for the calculation of successive
positions of a shock front is to develop a theory in which we can
find the shock strength, position and geometry of the shock at any
time without calculating the solution behind the shock. This is a
difficult task, because the nonlinear waves which
follow\footnote{Motion of a shock is influenced by the nonlinear waves which interact with the shock from both sides, ahead and behind it, but we consider here only the case when shock moves into a known
  state.} the shock strongly influence its evolution. Whitham
\cite{whitham-book} developed a simple and wonderful approximate
method, called the geometrical shock dynamics, in which the
effect of the flow behind the shock did not play any role. However, it
is well known that the fluid flow behind the shock has an important
effect on the shock motion and changes its strength, and this effect
has to be correctly accounted. Grinfel'd \cite{grinfeld} and Maslov
\cite{maslov} independently showed that the effect on the shock due to
the solution behind it can be described in the form of an infinite
system of transport equations for the normal derivatives of various
orders of a state variable behind the shock. This infinite system of
equations is obtained along the shock rays; see \cite{prasad-acta} for
more details of shock rays. Nevertheless, these transport equations
turn out to be quite involved and even for a weak shock, when the
equations simplify considerably, it is not easy to deal with an
infinite system of equations. As a remedy, Prasad and Ravindran
\cite{ravindran-prasad} proposed a procedure to truncate the infinite
system and thus developed the new theory of shock dynamics (NTSD). The
governing equations in this theory form a coupled set of differential
equations, consisting of the ray equations and those obtained by
truncating the infinite system. It is quite simple to use the NTSD to
one-dimensional (1-D) shock propagation, but its application to
multidimensional problems remained a challenge due to the formation of
kink (explained in the next two paragraphs) type of singularities,
which are found very frequently on a shock front. The appearance of
kinks were observed also in the numerical simulations of a nonlinear
wavefront\footnote{A nonlinear wavefront can be clearly distinguished
  from a shock front; see also \cite {arun-ct} for a more detailed
  explanation. For a comprehensive treatment and for references to the
  literature on the subject of this paper we refer the reader to
  \cite{prasad-book}.}; see \cite{ramanathan} for more details. The
kinks are formed on these fronts due to self-focusing. The first
experimental results showing the formation of kinks on a shock front
are available in the work of Sturtevant and Kulkarny
\cite{sturtevant-kulkarny}. These authors reported that the kinks
appear on a focusing shock front which is not too weak and it remained
a challenge to reproduce their experimental results by an analytical
method. However, Kevlahan \cite{kevlahan} used the NTSD to study this
problem and his numerical results agreed well not only with the
experimental results but also with some known exact and numerical
solutions of the Euler equations. This clearly shows the efficacy
of NTSD to produce physically realistic results.

A kink (kink line) is a point (curve) on a moving curve (surface),
across which the normal direction to the curve (surface) suffers a
jump discontinuity. The formation and propagation of kinks on a
shock front and also a nonlinear wavefront necessitated a
conservative formulation of the equations of NTSD as well as the
equations governing the motion of a weakly nonlinear wavefront. This
requirement finally led to the development of the kinematical
conservation laws (KCL) for the evolution of a curve in a plane by
Morton et al.\ \cite{morton} and for a surface in three-dimensional
(3-D) space by Giles et al.\ \cite{giles}.

For the sake of completeness and to make this paper
self-contained, in the following we briefly review the some basic
results from \cite{arun-ct,arun-prasad-wave}. The $d$-dimensional KCL
is a system of conservation equations governing the evolution of a
surface $\Omega_t$ in $\mathbb{R}^d$. The KCL is derived in
specially defined ray coordinates
$(\xi_1,\xi_2,\ldots,\xi_{d-1},t)$, where
$\xi_1,\xi_2,\ldots,\xi_{d-1}$ are the surface coordinates on
$\Omega_t$ and $t$ is time. The mapping between the ray
coordinates $(\xi_1,\xi_2,\ldots,\xi_{d-1},t)$ and the spatial
coordinates $(x_1,x_2,\ldots,x_d)$ is assumed to be locally
one-to-one. Since the KCL is a system of conservation laws, its
solutions may contain shocks in the ray coordinates. An image of any
one of these shocks, when mapped onto the
$(x_1,x_2,\ldots,x_d)$-space, is a kink surface on $\Omega_t$ across
which the normal direction to $\Omega_t$ and normal velocity are
discontinuous. Hence, the KCL is ideally suited to study the evolution
of a surface having kink type of singularities. However, the KCL being
purely a geometric result\footnote{KCL accounts for formation and propagation of kinks on a moving surface and the two compatibility conditions approximately incorporate results of Euler equations, see point 1. in appendix B}, it forms an incomplete system of equations
and additional closure relations are necessary to get a completely
determined set of equations. When $\Omega_t$ is a weakly nonlinear
wavefront in a polytropic gas \cite{prasad-1975}, the KCL system is
closed by a single equation representing the conservation of total
energy in a ray tube; see \cite{arun-etal-siam,arun-prasad-wave} for
more details. Throughout this paper we shall refer to the resulting
complete system of conservation laws, governing the evolution of a
nonlinear wavefront, as KCL based weakly nonlinear ray theory (WNLRT),
or briefly 3-D WNLRT in later sections. In the 2-D case this system
consists of just three conservation laws in $(\xi_1,t)$ coordinates,
which is hyperbolic when an appropriately defined non-dimensional
front velocity $m>1$; see \cite{prasad-book}. In the case of a
polytropic gas, $m>1$ corresponds to a wavefront on which the pressure
is greater than the constant pressure in the ambient gas. However, the
simplicity of 2-D KCL is lost when we consider the 3-D KCL, which is a
system of six conservation laws in $(\xi_1,\xi_2,t)$ coordinates with
three stationary divergence-free constraints. Following
\cite{arun-ct}, we shall refer to these three constraints together
as ``geometric solenoidal constraint''. The KCL based 3-D WNLRT
consists of seven equations and unlike the 2-D case, this system is
only weakly hyperbolic for $m>1$, in the sense that it has two
distinct eigenvalues and an eigenvalue of multiplicity five with a
four-dimensional eigenspace \cite{arun-prasad-wave,arun-prasad-eig}.

The eigenvalues of KCL based 3-D WNLRT are all real when
$m>1$. However, the weakly hyperbolic nature of the system and the
presence of geometric solenoidal constraint pose a challenge to
develop a numerical approximation. It is well known from the
literature that when the dimension of the eigenspace corresponding to
a multiple eigenvalue has a deficiency of one, the solution to a
Cauchy problem contains a mode, the so-called ``Jordan mode'', which
grows linearly in time. However, it has been shown in \cite{arun-ct}
that when the geometric solenoidal constraint is satisfied initially,
the solution to a Cauchy problem does not exhibit the Jordan mode.
Motivated by this, in \cite{arun-ct}, a constraint transport (CT)
technique has been built into a central finite volume scheme for the
3-D WNLRT system, therein the constraint is maintained up to machine
accuracy with the elimination of the Jordan mode. In addition, the
numerical method is shown to be very robust, second order accurate and
the numerical solution can be continued for a very long time, almost
indefinitely.

In \cite{baskar-jfm,monica} the authors have proposed a
conservative formulation of the NTSD based on 2-D KCL and it
turned out to be more effective than using the NTSD in a
differential form as done in \cite{kevlahan}. A conservative
formulation also has the advantage of using modern shock-capturing
algorithms and hence the kinks, whenever formed on the shock in
$(x_1,x_2)$-plane, can be automatically tracked as the images of
shocks in the ray coordinates $(\xi_1,t)$. The goal of the present
work is to derive the conservation laws of NTSD for a weak shock
based on the 3-D KCL and to use them for the computation of 3-D
shock fronts. We shall designate the conservation laws thus obtained
as the KCL based 3-D shock ray theory, or simply in 3-D SRT throughout
this paper. Our contribution consists in formulation and analysis of
3-D SRT, its numerical approximation and the results of numerical
experiments. These results reveal many realistic geometrical features
of shock fronts, the complex structure of kink lines and a comparison
with those in \cite{arun-ct} shows that a weak shock front and weakly
nonlinear wavefront are topologically the same.

The rest of this paper is organised as follows. In
section~\ref{sect:ntsd_eqns} we derive the system of conservation
laws of 3-D SRT and analyse its eigen-structure. In
section~\ref{sect:num_appxn} we briefly formulate a numerical
approximation of 3-D SRT, based on the constraint preserving, high
resolution central finite volume scheme from \cite{arun-ct}. The
results of numerical case studies are presented in
section~\ref{sect:case_study} and in order to facilitate the
comparison, the initial data for these tests are chosen as in
\cite{arun-ct}. In those test cases, where we observe the results
qualitatively similar to those in \cite{arun-ct}, we refrain from
presenting detailed and lengthy inferences. However, in a test
problem in which the interactions of kink lines and the corrugational
stability of a shock front are clearly observed, we present a
thorough analysis. In this case, as well as in the case of a radially
converging shock front, we highlight the important differences between
a weak shock front and a weakly nonlinear wavefront. Finally, we close
this paper with some concluding remarks in
section~\ref{sect:remarks}.

\section{Governing Equations of 3-D SRT}
\label{sect:ntsd_eqns}

A system of shock ray equations consists of the ray equations derived
from a shock manifold partial differential equation \cite{prasad-acta}
and an infinite system of compatibility conditions along a shock ray
\cite{anile-russo-1986,grinfeld,maslov,srinivasan-prasad-1985}. These
compatibility conditions are derived from the equations governing the
motion of the medium in which the shock propagates, e.g.\ the Euler
equations of gas dynamics. It has to be noted that unlike the well
known geometric optics theory for the propagation of a one-parameter
family of wavefronts, across which wave amplitudes are continuous, the
SRT with an infinite system of compatibility conditions is exact. This
is because the high frequency approximation required for the
derivation of jump relations is exactly satisfied for a shock
front.

Let us consider a shock propagating into a polytropic gas at rest and
in a uniform state $(\varrho,\uu{q},p)=(\varrho_0,\uu{0},p_0)$, where
$\varrho$ is the density, $\uu{q}=(q_1,q_2,q_3)$ is the particle
velocity and $p$ is the pressure. Let $a$ be the sound velocity in the
medium defined by $a^2=\gamma p/\varrho$, where $\gamma$ is the
ratio of specific heats. Let $\uu{N}$ denotes the unit normal to the
shock. Assuming the shock to be weak, the small amplitude
perturbations in the density $\varrho$, fluid velocity $\uu{q}$ and
pressure $p$ up to a short distance behind the shock can be expressed
using a small parameter $\varepsilon$ via the relations
\cite{prasad-book}
\begin{equation}
  \varrho-\varrho_0=\varepsilon\frac{\varrho_0}{a_0}\tilde{w}, \
  \uu{q}=\varepsilon \tilde{w}\uu{N}, \
  p-p_0=\varepsilon\varrho_0a_0\tilde{w},
\end{equation}
where $a_0$ is the sound speed in the uniform state and $\tilde{w}$ is
an amplitude of the order of unity, having the dimension of
velocity. Let us introduce a non-dimensional amplitude $\mu$, defined
on the shock front, via
\begin{equation}
  \mu:=\left.\frac{\tilde{w}}{a_0}\right\vert_{s}.
  \label{mu_defn}
\end{equation}
Following \cite{monica,prasad-book}, it can easily be seen that a
point $\uu{X}=(X_1,X_2,X_3)$ on a shock ray satisfies
\begin{align}
  \frac{\dd\uu{X}}{\dd T}&=a_0\left(1+\varepsilon\frac{\gamma+
    1}{4}\mu\right)\uu{N},\label{3d_srt_X}\\
  \frac{\dd\uu{N}}{\dd
    T}&=-\varepsilon\frac{\gamma+1}{4}a_0\uu{L}\mu,
  \label{3d_srt_N}
\end{align}
where $\dd/\dd T$ denotes the derivative
\begin{equation}
  \frac{\dd}{\dd T}:=\frac{\D}{\Dt}+a_0\left(1+\varepsilon\frac{\gamma+
    1}{4}\mu\right)\langle\uu{N},\nabla\rangle
  \label{d_dT_defn}
\end{equation}
and $\uu{L}$ is a tangential derivative along the shock front, defined
by
\begin{equation}
  \uu{L}:=\nabla-\uu{N}\langle\uu{N},\nabla\rangle.
  \label{3d_srt_L}
\end{equation}

For a weak shock, the first two of the infinite system of
compatibility conditions along the shock rays
(\ref{3d_srt_X})-(\ref{3d_srt_N}) are given by \cite{prasad-book}
\begin{align}
  \frac{\dd\mu}{\dd T}&=a_0{\it\Omega}_s\mu-\frac{\gamma+1}{4}\mu\mu_1,
  \label{3d_srt_mu1}\\
  \frac{\dd\mu_1}{\dd
    T}&=a_0{\it\Omega}_s\mu_1-\frac{\gamma+1}{2}\mu_1^2
  -\frac{\gamma+1}{4}\mu\mu_2.  \label{3d_srt_mu2}
\end{align}
At this point, we caution the reader that the coefficient of
$\mu_1^2$ in \eqref{3d_srt_mu2} has a misprint in the references
\cite{monica,prasad-book}. The term ${\it\Omega}_s$ in
\eqref{3d_srt_mu1}-\eqref{3d_srt_mu2}  denotes the mean curvature of
the shock and the variables $\mu_1$ and $\mu_2$ are defined by
\begin{equation}
  \mu_1:=\varepsilon\left.\langle\uu{N},\nabla\rangle
    \tilde{w}\right\vert_{s}, \
  \mu_2:=\varepsilon^2\left.\langle\uu{N},\nabla\rangle^2
    \tilde{w}\right\vert_{s}.
  \label{mu1_mu2}
\end{equation}
Due to the short wave approximation, both the quantities $\mu_1$ and
$\mu_2$ are of order unity \cite{monica,prasad-book}. The effect of
the term $\mu_1$ is very important in shock propagation. It represents
the gradient in the normal direction of the pressure or density just
behind the shock and takes into account of the effect of interactions
of the nonlinear waves which catch the shock from behind. It
has to be noted that apart from \eqref{3d_srt_mu1}-\eqref{3d_srt_mu2},
there exists an infinite system of compatibility conditions for
properly defined higher order derivatives $\mu_2,\mu_3,\mu_4,\ldots$ on the
shock. However, it appears to be very difficult to use this infinite
system of coupled equations for computing shock propagation.

The idea behind NTSD is to drop the term containing $\mu_2$ in
\eqref{3d_srt_mu2}; see \cite{ravindran-prasad}. Once this is done,
the system of equations \eqref{3d_srt_X}-\eqref{3d_srt_N} and
\eqref{3d_srt_mu1}-\eqref{3d_srt_mu2} is closed and it can then
used to compute shock propagation. For a discussion on the
validity of NTSD and its application to 2-D problems we refer the
reader to \cite{baskar-jfm,kevlahan,monica,prasad-book}. It is
interesting to note that the NTSD gives quite good results even
for a shock of arbitrary strength, which has been verified for a
1-D piston problem in \cite{lazarev-prasad-sing}. In addition, it
has also been observed in \cite{lazarev-prasad-sing} that the NTSD
takes less than $0.5\%$ of the computational time needed by a
typical finite difference method applied to the Euler equations\footnote{For some comments on the validity and accuracy of the NTSD, see point 2. in appendix B}.

\subsection{Conservation Forms of the Governing Equations}
As a first step, we non-dimensionalise all the independent and
dependent variables with the help of a characteristic length $L$
and the sound velocity $a_0$ in the uniform medium ahead of
the shock. We continue to denote all the resulting
non-dimensional variables also by the same symbols. Here, we
choose $L$ to be of the order of the distance over which the SRT
is valid and also the distance over which the shock propagates\footnote{For an explanation see point 3. in appendix B}. It
has to be noted that the nonlinear theory of Choquet-Bruhat
\cite{choquet-bruhat} is valid over a distance smaller than the
smaller of the radii of curvatures of the front at $t=0$. However,
for WNLRT and SRT, the distance $L$ could be far beyond the
caustic region as evident from the numerical results presented in
\cite{arun-etal-siam,baskar-jfm,monica,sangeeta}.

We now proceed to derive a conservation form the equations
\eqref{3d_srt_X}-\eqref{3d_srt_N} and
\eqref{3d_srt_mu1}-\eqref{3d_srt_mu2}. Let us define two variables $M$
and $\mathcal{V}$ via
\begin{equation}
  M:=1+\varepsilon\frac{\gamma+1}{4}\mu, \
  \mathcal{V}:=\frac{\gamma+1}{4}\mu_1.
  \label{M_V}
\end{equation}
Here, $M$ is a non-dimensional Mach number of the shock front and
$\mathcal{V}$ is an appropriately scaled normal derivative of the gas
density or pressure, just behind the shock. Following NTSD, we drop the
last term containing $\mu_2$ in (\ref{3d_srt_mu2}) and rewrite
(\ref{3d_srt_mu1})-(\ref{3d_srt_mu2}) in terms of $M$ and
$\mathcal{V}$ to obtain
\begin{align}
  \frac{\dd M}{\dd T}&={\it
    \Omega}_s(M-1)-\mathcal{V}(M-1),\label{3d_dMdT}\\
  \frac{\dd \mathcal{V}}{\dd T}&={\it
    \Omega}_s\mathcal{V}-2\mathcal{V}^2.\label{3d_dVdT}
\end{align}

We introduce a ray coordinate system $(\xi_1,\xi_2,t)$ on the
shock front in such a way that $t=\const$ is the shock front
$\Omega_t$ and ($\xi_1=\const,\xi_2=\const$) is a two-parameter
family of rays in the $(x_1,x_2,x_3)$-space. Let $\uu{U}$ and $\uu{V}$
be respectively the unit tangent vectors along the curves
($\xi_2=\const, t=\const$) and ($\xi_1=\const, t=\const$) on
$\Omega_t$. We denote by $G_1$ and $G_2$ respectively, the associated
metrics. On a given shock front at any time $t$ we have the relations
\begin{equation}
  \uu{X}_{\xi_1}=G_1\uu{U}, \ \uu{X}_{\xi_2}=G_2\uu{V},
  \label{X_xi1_X_xi2}
\end{equation}
where $G_1$ and $G_2$ are given by $G_1=\norm{\uu{X}_{\xi_1}}$ and
$G_2=\norm{\uu{X}_{\xi_2}}$. In view of the definition of $M$ from
(\ref{M_V}), in non-dimensional coordinates, the derivative $\dd/\dd
T$ defined in \eqref{d_dT_defn} assumes the form
\begin{equation}
  \label{eq:d_dT_M}
  \frac{\dd}{\dd T}=\frac{\D}{\Dt}+M\langle\uu{N},\nabla\rangle
\end{equation}
and it represents the time rate of change along a shock ray.
Hence, in the ray coordinates $(\xi_1,\xi_2,t)$, the derivative
$\dd/\dd T$ simply becomes $\D/\D t$.

As in \cite{arun-prasad-wave} it can be shown that the first and
second part of the ray equations, namely
(\ref{3d_srt_X})-(\ref{3d_srt_N}), are equivalent to the 3-D KCL
\begin{align}
  (G_1\uu{U})_t-(M\uu{N})_{\xi_1}&=0,\label{3d_srt_cons1}\\
  (G_2\uu{V})_t-(M\uu{N})_{\xi_2}&=0\label{3d_srt_cons2}
\end{align}
with the constraint
\begin{equation}
  (G_2\uu{V})_{\xi_1}-(G_1\uu{U})_{\xi_2}=0.
  \label{gsc}
\end{equation}
Therefore, it only remains to derive a conservation form of the
transport equations (\ref{3d_dMdT})-(\ref{3d_dVdT}). Our approach
is along the lines of \cite{prasad-parker} which follow a general
pattern valid for all compatibility conditions; see also
\cite{baskar-jfm}.

The ray tube area $\mathcal{A}$ of a tube of shock rays
\cite{prasad-book,whitham-book} is related to the mean curvature
${\it\Omega}_s$ by the relation
\begin{equation}
  \frac{1}{\mathcal{A}}\frac{\dd\mathcal{A}}{\dd l}=-2{\it \Omega}_s,
  \label{3d_srt_dAdl}
\end{equation}
where $l$ denotes the length measured along a shock ray. In
non-dimensional variables we have $\dd l=M\dd T$ and therefore from
(\ref{3d_srt_dAdl}) we get
\begin{equation}
  {\it \Omega}_s=-\frac{1}{2M\mathcal{A}}\frac{\dd\mathcal{A}}{\dd T}.
  \label{3d_srt_omegas}
\end{equation}
We notice the presence of the term ${\it\Omega}_s$ in both the
transport equations (\ref{3d_dMdT}) and (\ref{3d_dVdT}). In the light
of (\ref{3d_srt_omegas}), we conclude that ${\it \Omega}_s$ is related to
the rate of change of the ray tube area,  viz. $\dd\mathcal{A}/\dd
T$. Therefore, this term in the equations
(\ref{3d_dMdT})-(\ref{3d_dVdT}) represents the geometric decay or
amplification of the quantities $M$ and $\mathcal{V}$ .

We use the relation (\ref{3d_srt_omegas}) in (\ref{3d_dMdT})
with $\dd /\dd T$ replaced by $\D/\Dt$ in the ray coordinates and obtain
\begin{equation}
  \frac{2M}{M-1}M_t+\frac{\mathcal{A}_t}{\mathcal{A}}+2M\mathcal{V}=0.
  \label{3d_srt_MMt}
\end{equation}
Note that the left hand side of the above equation gives a combination
$\{f^\prime(h)/f(h)\}h_t+\mathcal{A}_t/\mathcal{A}$, where $h=M-1$ and
$f$ is the function $f(h):=h^2e^{2h}$. Hence, we get a conservation form
\begin{equation}
  \left\{\mathcal{A}f(M-1)\right\}_t+2\mathcal{A}Mf(M-1)\mathcal{V}=0.
  \label{3d_srt_consF}
\end{equation}
The ray tube area $\mathcal{A}$ is given by the expression
\begin{equation}
  \mathcal{A}=G_1G_2\sin\Psi,
\end{equation}
where $\Psi$ is the angle between the two unit vectors $\uu{U}$ and
$\uu{V}$. Therefore, from (\ref{3d_srt_consF}), we finally get a
balance equation
\begin{equation}
  \left\{(M-1)^2e^{2(M-1)}G_1G_2\sin\Psi\right\}_t
  +2M(M-1)^2e^{2(M-1)}G_1G_2\mathcal{V}\sin\Psi=0.
  \label{3d_srt_cons3}
\end{equation}

Similarly, using the expression (\ref{3d_srt_omegas}) for ${\it
  \Omega}_s$ in (\ref{3d_dVdT}) we obtain an equation which we rewrite
as
\begin{equation}
  \mathcal{V}_t+\frac{\mathcal{V}}{2\mathcal{A}}\mathcal{A}_t
  +\frac{\mathcal{V}}{2}\left(\frac{1}{M}-1\right)
  \frac{\mathcal{A}_t}{\mathcal{A}}+2\mathcal{V}^2=0.
  \label{3d_srt_A}
\end{equation}
We use (\ref{3d_srt_MMt}) to replace the factor
$\mathcal{A}_t/\mathcal{A}$ from the third term in (\ref{3d_srt_A})
and write the resulting equation as
\begin{equation}
  \left\{\ln\left(\mathcal{V}^2\mathcal{A}\right)+2(M-1)\right\}_t
  +(M+1)\mathcal{V}=0,
\end{equation}
which gives a balance equation
\begin{equation}
  \left\{e^{2(M-1)}G_1G_2\mathcal{V}^2\sin\Psi\right\}_t
  +(M+1)e^{2(M-1)}G_1G_2\mathcal{V}^3\sin\Psi=0.
  \label{3d_srt_cons4}
\end{equation}
\begin{remark}
  \label{rem:terms}
  The second term in \eqref{3d_srt_MMt} and the second and third terms
  in \eqref{3d_srt_A} represent the geometrical effect of convergence
  or divergence of rays. The third term in \eqref{3d_srt_MMt}
  represents the effect of interaction of nonlinear waves which
  overtake the shock from behind. The fourth term in \eqref{3d_srt_A}
  is the usual effect of genuine nonlinearity which governs the
  evolution of $\mathcal{V}$, also seen in the 1-D model
  $u_t+uu_x=0$.
\end{remark}
The complete set of equations of 3-D KCL based NTSD,
hereafter designated as the conservation laws of 3-D SRT, consists of
the equations (\ref{3d_srt_cons1})-(\ref{3d_srt_cons2}),
(\ref{3d_srt_cons3}) and (\ref{3d_srt_cons4}).
\begin{remark}
  The conservation forms the \eqref{3d_srt_cons3} and \eqref{3d_srt_cons4} of the compatibility conditions (\ref{3d_dMdT}) and (\ref {3d_dVdT}) respectively, are physically realistic. They are the three-dimensional extensions of those
  derived by Baskar and Prasad \cite{baskar-jfm}, namely
  \begin{align}
    \left\{(M-1)^2e^{2(M-1)}G\right\}_t+2M(M-1)^2e^{2(M-1)}G\mathcal{V}&=0,
    \label{2d_srt_cons3}\\
     \left\{e^{2(M-1)}G\mathcal{V}^2\right\}_t+(M+1)e^{2(M-1)}G\mathcal{V}^3&=0.
     \label{2d_srt_cons4}
  \end{align}
  Note that \eqref{2d_srt_cons3} and \eqref{2d_srt_cons4} follow from
  \eqref{3d_srt_cons3} and \eqref{3d_srt_cons4} respectively, by
  replacing $G_1G_2\sin\Psi$ by $G$. In the linear theory, the energy
  conservation along a ray tube\footnote{Ray tube has been sketched and explained in section 3 of [3].} is represented by $\{G(M-1)^2\}_t=0$,
  which can be written in an integral formulation using two
  cross-sections of a ray tube; see \cite{whitham-book}. The genuine nonlinearity in Euler equations stretches a shock ray (i.e., makes the shock front move faster) due to presence of the term $M\langle\uu{N},\nabla\rangle$ in (2.14) and hence the flux of energy across a section of a shock ray tube increases by a factor $e^{2(M-1)}$. Though this term looks small for a weak shock, its accumulative effect over long time is significant. Presence of  this term is also seen in the conservation laws of 3-D WNLRT in \cite{arun-prasad-wave}. The
  source term in \eqref{3d_srt_cons3}, leading to a decay in the shock
  velocity when $\mathcal{V}>0$, is the main cause of differences in
  the results as compared to the results for a nonlinear wavefront. We
  shall postpone a detailed discussion on this to
  section~\ref{sect:case_study}.
\end{remark}

One of the aims of this paper is to compare the the results of 3-D SRT
with those of 3-D WNLRT, reported in \cite{arun-ct}. Hence, for the
sake of completeness, in the following we reproduce the conservation
laws 3-D WNLRT, i.e.\
\begin{align}
  (g_1\uu{u})_t-(m\uu{n})_{\xi_1}&=0, \label{3d_wnlrt_cons1}\\
  (g_2\uu{v})_t-(m\uu{n})_{\xi_2}&=0, \label{3d_wnlrt_cons2}\\
  \left((m-1)^2e^{2(m-1)}g_1g_2\sin\psi\right)_t&=0 \label{3d_wnlrt_cons3}
\end{align}
with the constraint
\begin{equation}
  \label{3d_wnlrt_cons4}
  (g_2\uu{v})_{\xi_1}-(g_1\uu{u})_{\xi_2}=0.
\end{equation}
We notice that unlike the balance equations of 3-D SRT, the system of
conservation laws \eqref{3d_wnlrt_cons1}-\eqref{3d_wnlrt_cons3} is
homogeneous, i.e.\ without any source term.

\subsection{Eigen-structure of the System of Conservation Laws}
This section is devoted to the analysis of the eigenvalues and
eigenvectors of the system conservations laws
(\ref{3d_srt_cons1})-(\ref{3d_srt_cons2}), (\ref{3d_srt_cons3}) and
(\ref{3d_srt_cons4}). First, we recast them in the usual divergence form
\begin{equation}
  W_t+{F_1(W)}_{\xi_1}+{F_2(W)}_{\xi_2}=S(W)
  \label{3d_srt_cons}
\end{equation}
with the conserved variable $W$, the fluxes $F_1(W),F_2(W)$ and
the source term $S(W)$, given as
\begin{equation}
  \begin{aligned}
    W&=\left(G_1\uu{U},G_2\uu{V},(M-1)^2e^{2(M-1)}G_1G_2\sin\Psi,
    e^{2(M-1)}G_1G_2\mathcal{V}^2\sin\Psi\right)^T,\\
    F_1(W)&=\left(M\uu{N},\uu{0}, 0, 0\right)^T,\\
    F_2(W)&=\left(\uu{0}, M\uu{N}, 0, 0\right)^T,\\
    S(W)&=\left(\uu{0},\uu{0},-2M(M-1)^2e^{2(M-1)}G_1G_2\mathcal{V}\sin\Psi,
    -(M+1)e^{2(M-1)}G_1G_2\mathcal{V}^3\sin\Psi\right)^{T}.
  \end{aligned}
  \label{W_F1_F2_S}
\end{equation}
In order to discuss the eigenvalues and eigen-structure of
\eqref{3d_srt_cons}, we need its explicit form as a system of partial
differential equations. Introducing a vector of primitive variables
via $V=(U_1,U_2,V_1,V_2,M,G_1,G_2,\mathcal{V})^T$, we can derive a
quasi-linear form of \eqref{3d_srt_cons} as
\begin{equation}
  \tilde{A}V_t+\tilde{B}^{(1)}V_{\xi_1}+\tilde{B}^{(2)}V_{\xi_2}=\tilde{C},
  \label{3d_srt_quasi}
\end{equation}
where the expressions for the matrices
$\tilde{A},\tilde{B}^{(1)},\tilde{B}^{(2)}$ and $\tilde{C}$ are given
in Appendix~\ref{sect:append_a}. It is interesting to note that the
matrices $\tilde{A},\tilde{B}^{(1)}$ and $\tilde{B}^{(2)}$ admit the
following block structure
\begin{equation}
  \tilde{A}=
  \begin{pmatrix}
    A & O_{7,1}\\
    R_{1,7} & \frac{2G_1G_2}{\mathcal{V}}
  \end{pmatrix}, \
  \tilde{B}^{(1)}=
  \begin{pmatrix}
    B^{(1)} & O_{7,1}\\
    O_{1,7} & 0
  \end{pmatrix}, \
  \tilde{B}^{(2)}=
  \begin{pmatrix}
    B^{(2)} & O_{7,1}\\
    O_{1,7} & 0
  \end{pmatrix}.
  \label{mat_block}
\end{equation}
Here, $A,B^{(1)}$ and $B^{(2)}$ are the corresponding flux Jacobian
matrices of 3-D WNLRT, cf.\ \cite{arun-etal-siam,arun-prasad-wave},
$O$ denotes a zero-matrix and $R_{1,7}$ is a row-matrix; see also
Appendix~\ref{sect:append_a}. The characteristic equation of
\eqref{3d_srt_quasi} is given by
\begin{equation}
  \det\tilde{M}_{8,8}(\lambda)\equiv\det\left(e_1\tilde{B}^{(1)}
    +e_2\tilde{B}^{(2)} -\lambda\tilde{A}\right)=0.
  \label{3d_srt_char_eqn}
\end{equation}
Using the block structure of the matrices in \eqref{mat_block} we can
obtain
\begin{equation}
  \tilde{M}_{8,8}(\lambda)=
  \begin{pmatrix}
    M_{7,7}(\lambda) & O_{7,1}\\
    -\lambda R_{1,7} & -\lambda\frac{2G_1G_2}{\mathcal{V}}
  \end{pmatrix},
  \label{m_88}
\end{equation}
where $M_{7,7}(\lambda):=e_1B^{(1)}+e_2B^{(2)}-\lambda A$ is the
matrix pencil of the 3-D WNLRT. Therefore, the characteristic equation
\eqref{3d_srt_char_eqn} simplifies as
\begin{equation}
  \lambda\det M_{7,7}(\lambda)=0.
\end{equation}
Using the results of \cite{arun-prasad-wave,arun-prasad-eig}, it is
easy to find the roots of the polynomial equation
$\det M_{7,7}(\lambda)=0$ in $\lambda$ and the nullvectors of
$M_{7,7}$ corresponding to these roots. The roots of
\eqref{3d_srt_char_eqn} can be obtained as
$\lambda_1,\lambda_2(=-\lambda_1),\lambda_3=\cdots=\lambda_8=0$, where
\begin{equation}
  \lambda_1=\left\{\frac{M-1}{2\sin^2\Psi}\left(\frac{e_1^2}{G_1^2}-
  \frac{2e_1e_2}{G_1G_2}\cos\Psi+\frac{e_2^2}{G_2^2}\right)\right\}^{\frac{1}{2}}
  \label{3d_srt_lamb1}
\end{equation}
and $(e_1,e_2)\in\mathbb{R}^2$ with $e_1^2+e_2^2=1$. Let us now
consider the matrix $\tilde{M}_{8,8}(0)$ for the multiple eigenvalue
$\lambda=0$. Clearly, the rank of the matrix $\tilde{M}_{8,8}(0)$ is
same as that of $M_{7,7}(0)$, cf.\ \eqref{m_88}. It has been shown in
\cite{arun-prasad-wave} that the rank of $M_{7,7}(0)$ is three. Thus,
the dimension of the nullspace of $\tilde{M}_{8,8}(0)$ is
$8-3=5$. This proves that the dimension of the eigenspace
corresponding to $\lambda=0$ of multiplicity six is only five. We
summarise these results as a theorem.
\begin{theorem}
  \label{thm:3d_srt_eig}
  The system (\ref{3d_srt_quasi}) has eight eigenvalues $\lambda_1,
  \lambda_2~(=-\lambda_1), \lambda_3=\lambda_4=\cdots=\lambda_8=0$,
  where $\lambda_1$ and $\lambda_2$ are real for $M>1$ and purely
  imaginary for $M<1$. Further, the dimension of the eigenspace
  corresponding to the multiple eigenvalue zero is five.
\end{theorem}
\begin{remark}
  It is important to note that $M<1$ is not physically unrealistic for a shock. For example, the signature of a sonic boom produced by a convex and smooth upper surface an aerofoil with a
sharp leading edge (see figures 2.1 and 2.2, \cite{baskar-sonic}) consists of a leading shock followed by a  continuous flow in which the pressure decreases and terminates in a trailing shock. In the forward part
    of the continuous flow the pressure is greater than that in the
    ambient medium and in the rear part it is less than that behind
    the trailing shock. The shock velocity of the trailing weak shock
    is half of the sound velocity behind it, say $a_0$ and the sound
    speed ahead of it, which is less than $a_0$. Thus Mach number $M$
    of the trailing shock is less than $1$. 
\end{remark}

\section{Numerical Approximation}
\label{sect:num_appxn}

As a consequence of Theorem~\ref {thm:3d_srt_eig} we infer that the
system of conservation laws \eqref{3d_srt_cons} is only weakly
hyperbolic for $M>1$; hence, an initial value problem may not be
well-posed in the strong hyperbolic sense. In addition, weakly
hyperbolic systems are likely to be more sensitive than regular
hyperbolic systems also from a computational point of view. Numerical
as well as theoretical analysis indicates that their solution may not
belong to the BV spaces and can only be measure valued. Despite such
theoretical difficulties, in \cite{arun-ct,arun-etal-siam}, we have
been able to develop accurate and efficient numerical approximations
of the analogous weakly hyperbolic system of 3-D WNLRT using, simple but
robust, central schemes. Since the conservation laws of 3-D SRT and
those of 3-D WNLRT, viz.\
\eqref{3d_wnlrt_cons1}-\eqref{3d_wnlrt_cons4}, are structurally
similar, we do not intend to present the details of the numerical
approximation and the refer the reader to \cite{arun-ct} for more
details.

We discretise the system \eqref{3d_srt_cons} using the cell integral
averages $\overline{W}_{i,j}$ of the conservative variable $W$, taken
over square mesh cells. From the given cell averages
$\overline{W}_{i,j}^n$ at time $t^n$, we reconstruct a piecewise
linear interpolant using the standard MUSCL type procedures. In order
to obtain the discrete slopes in the $\xi_1$- and $\xi_2$-directions, we
employ a central weighted essentially non-oscillatory limiter
\cite{jiang-shu}. The piecewise linear reconstruction enables us to
compute the cell interface values of the conserved variable $W$.

The starting point for the construction of numerical scheme is a
semi-discrete discretisation of (\ref{3d_srt_cons}), given by
\begin{equation}
  \frac{\dd\overline{W}_{i,j}}{\dd t}=-\frac{{\mathcal{F}_1}_{i+\frac{1}{2},j}
    -{\mathcal{F}_1}_{i-\frac{1}{2},j}}{h_1}
  -\frac{{\mathcal{F}_2}_{i,j+\frac{1}{2}}
    -{\mathcal{F}_2}_{i,j-\frac{1}{2}}}{h_2}+S(\overline{W}_{i,j}),
  \label{2d_semi_disc}
\end{equation}
where the quantities $\mathcal{F}_{i+1/2,j}$ and
$\mathcal{F}_{i,j+1/2}$ are respectively, the numerical fluxes at
the cell interfaces $(i+1/2,j)$ and $(i,j+1/2)$. We employ the high
resolution flux given by Kurganov and Tadmor \cite{kurganov-tadmor},
for these interface fluxes, e.g.\ at a right-hand vertical edge
\begin{equation}
  {\mathcal{F}_1}_{i+\frac{1}{2},j}\left(W_{i,j}^R,W_{i+1,j}^L\right)
  =\frac{1}{2}\left(F_1\left(W_{i+1,j}^L\right)+F_1\left(W_{i,j}^R\right)\right)
  -\frac{a_{i+\frac{1}{2},j}}{2}\left(W_{i+1,j}^L-W_{i,j}^R\right),
  \label{2d_fluxes}
\end{equation}
where $W_{i,j}^{L(R)}$ denote respectively, the left and right
interpolated states at the interface $(i+1/2,j)$. The expression for
the flux $\mathcal{F}_{i,j+1/2}$ at an upper horizontal edge is
analogous. In the flux formula \eqref{2d_fluxes}, the term
$a_{i+1/2,j}$ denotes the local speed of propagation at cell
interfaces, given by
\begin{equation}
  a_{i+\frac{1}{2},j}:=\max\left\{\rho\left(\frac{\D F_1}{\D
        W}\left(W_{i,j}^R\right)\right),\rho\left(\frac{\D F_1}{\D
        W}\left(W_{i+1,j}^L\right)\right)\right\}
  \label{2d_local_speed}
\end{equation}
where $\rho(A):=\max_i\abs{\lambda_i(A)}$, with $\lambda_i(A)$ being the
eigenvalues of the matrix $A$.

To improve the temporal accuracy and to gain second order accuracy in
time, we use a TVD Runge-Kutta scheme \cite{shu-tvd-rk} to numerically
integrate the system of ordinary differential equations in
(\ref{2d_semi_disc}). Denoting the right hand side of
(\ref{2d_semi_disc}) by $\mathcal{L}_{i,j}(W)$, the second order
Runge-Kutta scheme updates $W$ through the following two stages
\begin{equation}
  \begin{aligned}
    W_{i,j}^{(1)}&=\overline{W}_{i,j}^n+\dt\mathcal{L}_{i,j}
    \left(\overline{W}^n\right),\\
    \overline{W}_{i,j}^{n+1}&=\frac{1}{2}\overline{W}_{i,j}^n
    +\frac{1}{2}W_{i,j}^{(1)}+\frac{1}{2}\dt\mathcal{L}_{i,j}\left(W^{(1)}\right).
  \end{aligned}
  \label{2d_tvd_rk}
\end{equation}

It is to be noted that any consistent numerical solution of the 3-D
KCL system \eqref{3d_srt_cons1}-\eqref{3d_srt_cons2} also has to
satisfy the geometric solenoidal constraint \eqref{gsc} at any
time. Note that this constraint is an involution for the 3-D KCL
\eqref{3d_srt_cons1}-\eqref{3d_srt_cons2}, i.e.\ once fulfilled at the
initial data it is fulfilled for all times. Since the physically exact
solution has this feature, a numerical solution should also possess it some
discrete sense. Hence, in the numerical approximation of the analogous
3-D WNLRT system in \cite{arun-ct}, a CT algorithm \cite{evans-hawley}
was built into the central finite volume to enforce the geometric
solenoidal constraint. In what follows, we briefly review this CT
strategy and refer to \cite{arun-ct} for more details on its
implementation.

The geometric solenoidal constraint \eqref{gsc} guarantees the
existence of three potential functions $\mathbb{A}_k,k=1,2,3$, such
that
\begin{equation}
  \label{eq:A_xi1_xi2}
  G_1U_k={\mathbb{A}_k}_{\xi_1}, \ G_2V_k={\mathbb{A}_k}_{\xi_2}.
\end{equation}
Using \eqref{eq:A_xi1_xi2} in the 3-D KCL system
\eqref{3d_srt_cons1}-\eqref{3d_srt_cons2} yields the evolution
equations
\begin{equation}
  \label{eq:A_evol}
  {\mathbb{A}_k}_t=MN_k.
\end{equation}
In the CT method, we store the three potentials $\mathbb{A}_k$ at the
centres of a staggered grid. With aid of these potentials we redefine
the values of the vectors $G_1\uu{U}$ and $G_2\uu{V}$ at the cell
edges, which are treated as length averaged quantities. In this
way, $G_2\uu{V}$ is collocated at $(i+1/2,j)$, whereas $G_1\uu{U}$
is collocated at $(i,j+1/2)$. The definitions of these collocated
values are obtained by simply discretising the derivatives in
\eqref{eq:A_xi1_xi2} using central differences, i.e.\
\begin{align}
  [G_1U_k]_{i,j+\frac{1}{2}}&=\frac{1}{h_1}\left({\mathbb{A}_k}
    _{i+\frac{1}{2},j+\frac{1}{2}}-{\mathbb{A}_k}_{i-\frac{1}{2},j+\frac{1}{2}}\right),
  \label{g1u_defn}\\
  [G_2V_k]_{i+\frac{1}{2},j}&=\frac{1}{h_2}\left({\mathbb{A}_k}
    _{i+\frac{1}{2},j+\frac{1}{2}}-{\mathbb{A}_k}_{i+\frac{1}{2},j-\frac{1}{2}}\right)
  \label{g2v_defn}.
\end{align}
With the above collocated values we calculate
\begin{equation}
  \label{eq:gsc_centre}
  \left.(G_2V_k)_{\xi_1}-(G_1U_k)_{\xi_2}\right\vert_{i,j}=\frac{1}{h_1}\left(
    [G_2V_k]_{i+\frac{1}{2},j}-[G_2V_k]_{i-\frac{1}{2},j}\right)-\frac{1}{h_2}\left(
    [G_1U_k]_{i,j+\frac{1}{2}}-[G_1U_k]_{i,j-\frac{1}{2}}\right).
\end{equation}
Using \eqref{g1u_defn}-\eqref{g2v_defn}, we can easily see that the
right hand side of \eqref{eq:gsc_centre} vanishes due to perfect
cancellation. Hence, in this way, we have devised a method to enforce
the geometric solenoidal constraint at the cell centres of the finite
volumes.

It remains to be specified how to compute the values of the potentials
$\mathbb{A}_k$ on the staggered grids. Note that integrating
\eqref{eq:A_evol} over a staggered grid yields the update formula
\begin{equation}
  \label{eq:a_update}
  \frac{\dd}{\dd
    t}{\mathbb{A}_k}_{i+\frac{1}{2},j+\frac{1}{2}}=[MN_k]_{i+\frac{1}{2},j+\frac{1}{2}}.
\end{equation}
In order to compute the expression on the right hand side we use a
simple averaging, i.e.\
\begin{equation}
  \label{eq:a_avg}
  [MN_k]_{i+\frac{1}{2},j+\frac{1}{2}}=\frac{1}{4}\left([MN_k]_{i+\frac{1}{2},j}
    +[MN_k]_{i,j+\frac{1}{2}}+[MN_k]_{i+\frac{1}{2},j+1}
    +[MN_k]_{i+1,j+\frac{1}{2}}\right),
\end{equation}
where the values of $MN_k$ on the cell edges are obtained from the
numerical fluxes $\mathcal{F}_1$ and $\mathcal{F}_2$ of the finite
volume scheme, cf. also \eqref{W_F1_F2_S}. The resulting ordinary
differential equations in \eqref{eq:a_update} are integrated using the
same TVD Runge-Kutta method in \eqref{2d_tvd_rk}. At the beginning of
the next time step, the cell centred values of $G_1U_k$ and $G_2V_k$
are calculated by interpolation
\begin{align}
  [G_1U_k]_{i,j}&=\frac{1}{2}\left([G_1U_k]_{i,j+\frac{1}{2}}
    +[G_1U_k]_{i,j-\frac{1}{2}}\right), \label{eq:g1u_interp}\\
  [G_2V_k]_{i,j}&=\frac{1}{2}\left([G_2V_k]_{i+\frac{1}{2},j}
    +[G_2V_k]_{i-\frac{1}{2},j}\right). \label{eq:g2v_interp}
\end{align}

There is an extensive discussion in
\cite{arun-ct,arun-prasad-wave} on the formulation of the
initial data and the implementation of boundary conditions needed for
the numerical method. Since the variables $G_1,G_2,\uu{U},\uu{V}$ and
$M$ are analogous for 3-D SRT and 3-D WNLRT, we do not discuss the initial
conditions for these variables. The only additional unknown in 3-D
SRT, namely the gradient $\mathcal{V}$ of the gas density at the shock
front, has been assigned constant positive value
$\mathcal{V}=\mathcal{V}_0$.

At the end of each time step, we get the updated value of the
conserved variable $W$. Since $\norm{\uu{U}}=\norm{\uu{V}}=1$, from
the first six components of $W$, the values of $G_1,G_2,\uu{U}$ and
$\uu{V}$ can be computed very easily. To get the updated value of the
normal velocity $M$ we proceed as follows. Note that
\begin{equation}
  (M-1)^2e^{2(M-1)}=\frac{W_7}{G_1G_2\sin\Psi}\equiv \kappa, \ \mbox{say}.
  \label{M_eqn}
\end{equation}
We now solve for $M$ the nonlinear equation
\begin{equation}
  \vartheta(M)\equiv (M-1)^2e^{2(M-1)}-\kappa=0
  \label{theta}
\end{equation}
using Newton-Raphson method. The monotonicity of the function
$\vartheta$ in $(1,\infty)$ ensures the uniqueness of the solution of
(\ref{theta}). Having got the values of $G_1,G_2,\uu{U},\uu{V}$ and
$M$, the value of $\mathcal{V}$ can be easily obtained from $W_8$. In
order to get the successive positions of the shock front $\Omega_t$,
we numerically integrate the first part, viz.\ \eqref{3d_srt_X}, of
the ray equations, i.e.\
\begin{equation}
  \label{eq:ray_num}
  \frac{\dd }{\dd t}\uu{X}_{i,j}(t)=M_{i,j}(t)\uu{N}_{i,j}(t)
\end{equation}
using the two-stage Runge-Kutta method; see also
\cite{arun-ct,arun-etal-siam} for more details.

\section{Numerical Case Studies}
\label{sect:case_study}

In this section we present the evolution of shock fronts, obtained
by the results of extensive numerical simulations of the balance
laws \eqref{3d_srt_cons}, starting with a wide range of initial
geometries. Our motivation for these numerical experiments is
multifold.\ (i) Demonstration the efficiency of SRT to produce
intricate the shapes of shocks with complex patterns of kink lines
formed on them, (ii) Qualitative assessment of the numerical
results with the existing experimental and numerical results in
the literature, (iii) A comparison of the geometrical shapes of a weak
shock front and a weakly nonlinear wavefront. Our intention is to
bring out the important similarities and differences between a weak
shock front and weakly nonlinear wavefront and in order to facilitate
the comparison, we choose the same set of initial data as given in
\cite{arun-ct}, except for an additional test problem reported in
subsection \ref{subsec:axisym}. As remarked earlier in the
introduction, in those test problems where we observe qualitatively
similar geometries, we do not intent to present a detailed
discussion. However, we clearly point out those aspects where a shock
front differs significantly from a nonlinear wavefront in its
evolution.

In all the numerical case studies performed below, we have used the
CFL condition
\begin{equation}
  \nu=\dt\max\left(\frac{\rho_1}{h_1},\frac{\rho_2}{h_2}\right),
  \label{cfl_cond}
\end{equation}
where $\rho_1$ and $\rho_2$ are respectively the maxima of the
absolute values of the eigenvalues of the flux Jacobian matrices in
$\xi_1$- and $\xi_2$-directions, cf.\ \eqref{3d_srt_lamb1} and the
maxima are taken over all mesh cells. We have set $\nu=0.9$ in all our
numerical experiments. Further, in all the test problems, the normal
derivative of the density at the shock, viz.\ $\mathcal{V}$, has been
initialised with a constant positive value $\mathcal{V}_0=0.2$. This
is due to the fact the gradient of the pressure behind the shock is
positive for almost all physically realistic shock fronts; see also
\cite{monica} for a related discussion.
\begin{remark}
  In all our numerical experiments the CT technique preserves the
  geometric solenoidal constraint \eqref{gsc} up to machine accuracy
  for all times. In Figure~\ref{dip_srt_div} we plot the discrete
  divergence of $\mathfrak{B}_1:=(G_2V_1,-G_1U_1)$ at time
  $t=10.0$, for the problem studied in subsection
  \ref{subsec:non_axisym}, and the figure shows the error to be of
  order $10^{-15}$. Since in all other test cases we observe the same
  behaviour, we refrain from presenting the details. The discrete
  divergences of the remaining two components $\mathfrak{B}_2$ and
  $\mathfrak{B}_3$ also show the same trend.
\end{remark}

\subsection{ Propagation of a Non-axisymmetric Shock Front}
\label{subsec:non_axisym}

We choose the initial shock front $\Omega_0$ in a such a way that it
is not axisymmetric. The front $\Omega_0$ has a single smooth dip with
the initial shape given by
\begin{equation}
  \Omega_0\colon
  x_3=\frac{-\kappa}{1+\frac{x_1^2}{\alpha^2}+\frac{x_2^2}{\beta^2}},
  \label{dip_omega0}
\end{equation}
where the parameter values are set to
$\kappa=0.5,\alpha=1.5,\beta=3$.

\begin{figure}
  \centering
  \begin{tabular}{cc}
    \includegraphics[width=0.5\textwidth]{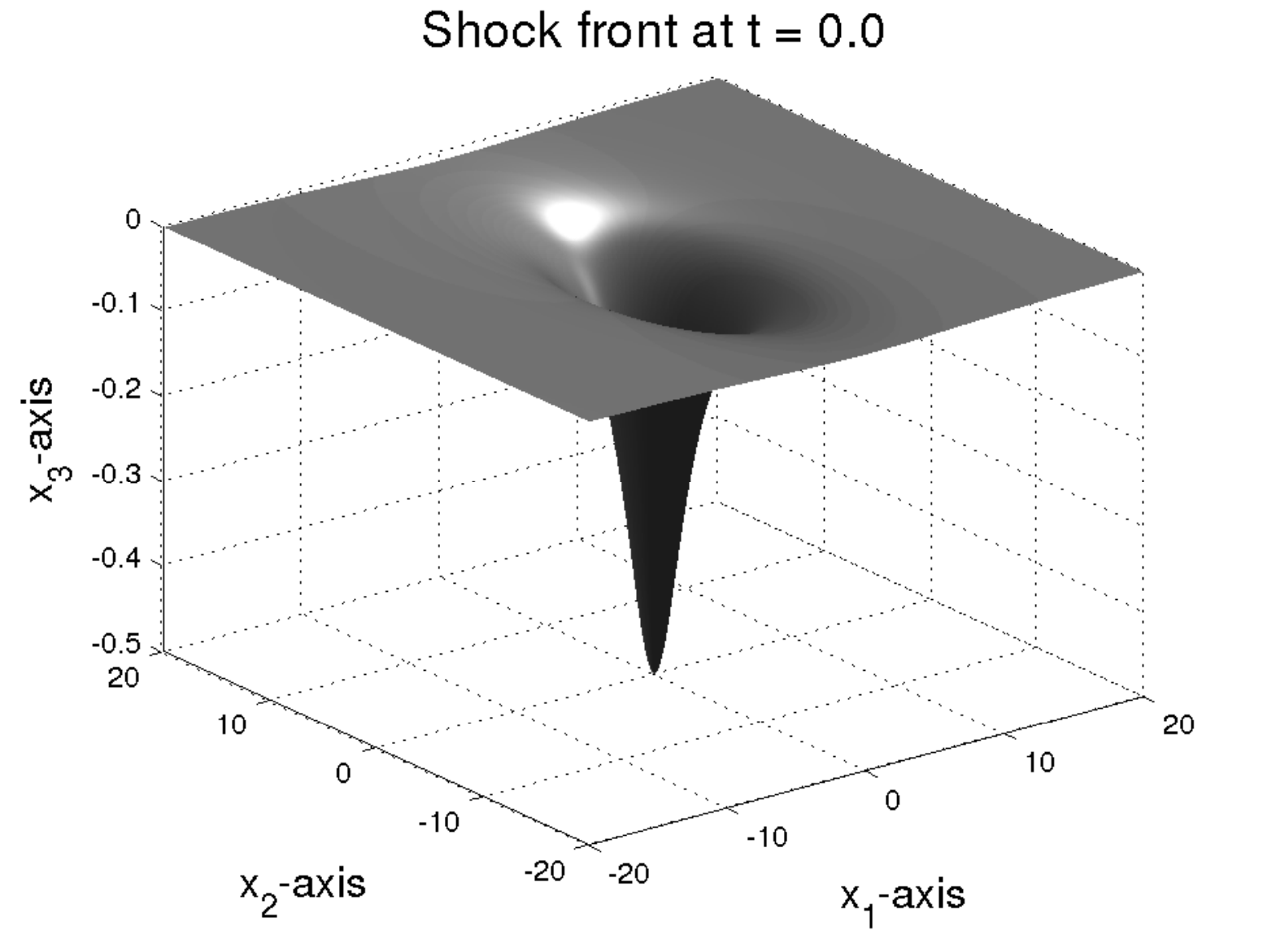}&
    \includegraphics[width=0.5\textwidth]{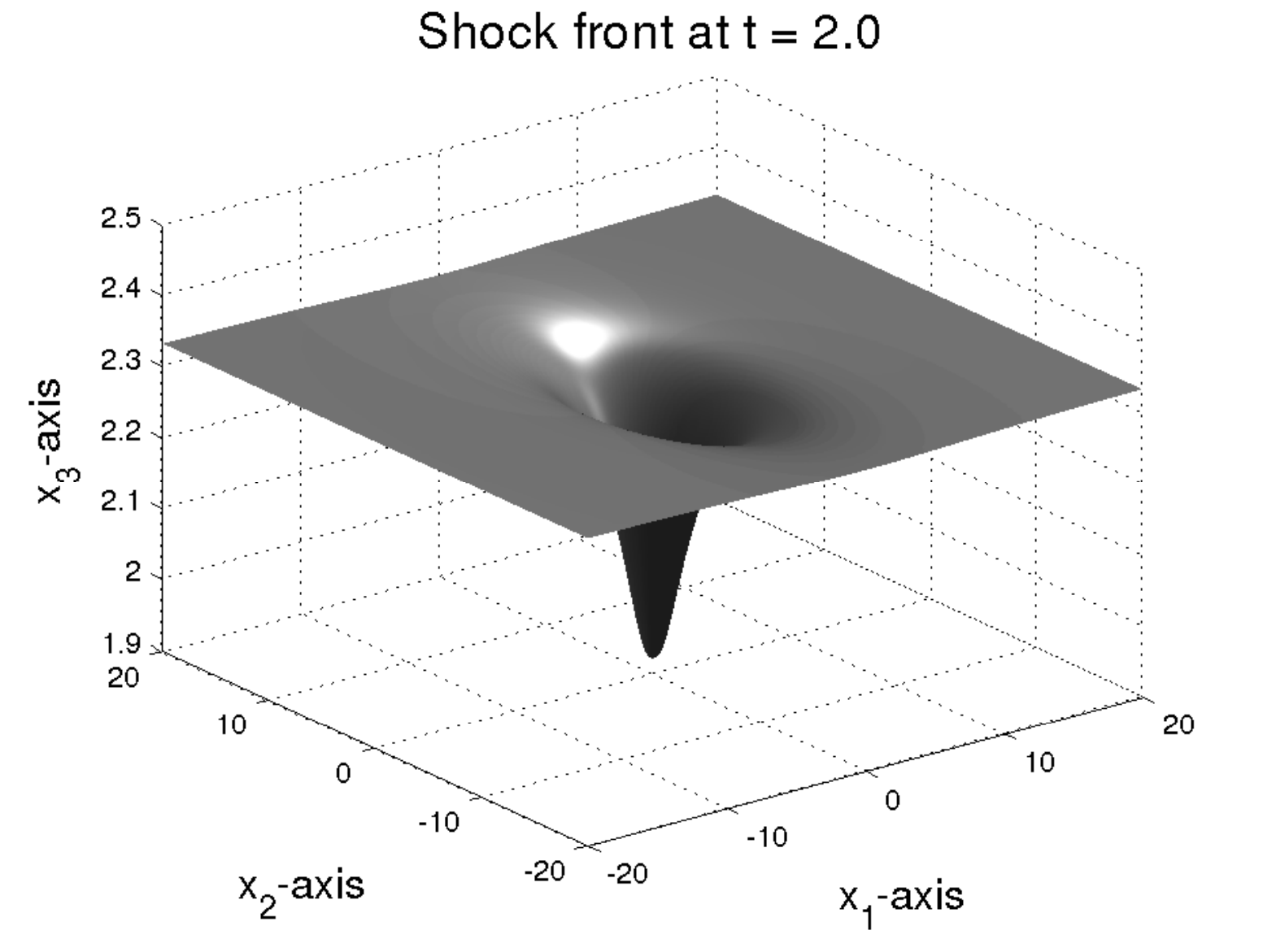}\\
    \includegraphics[width=0.5\textwidth]{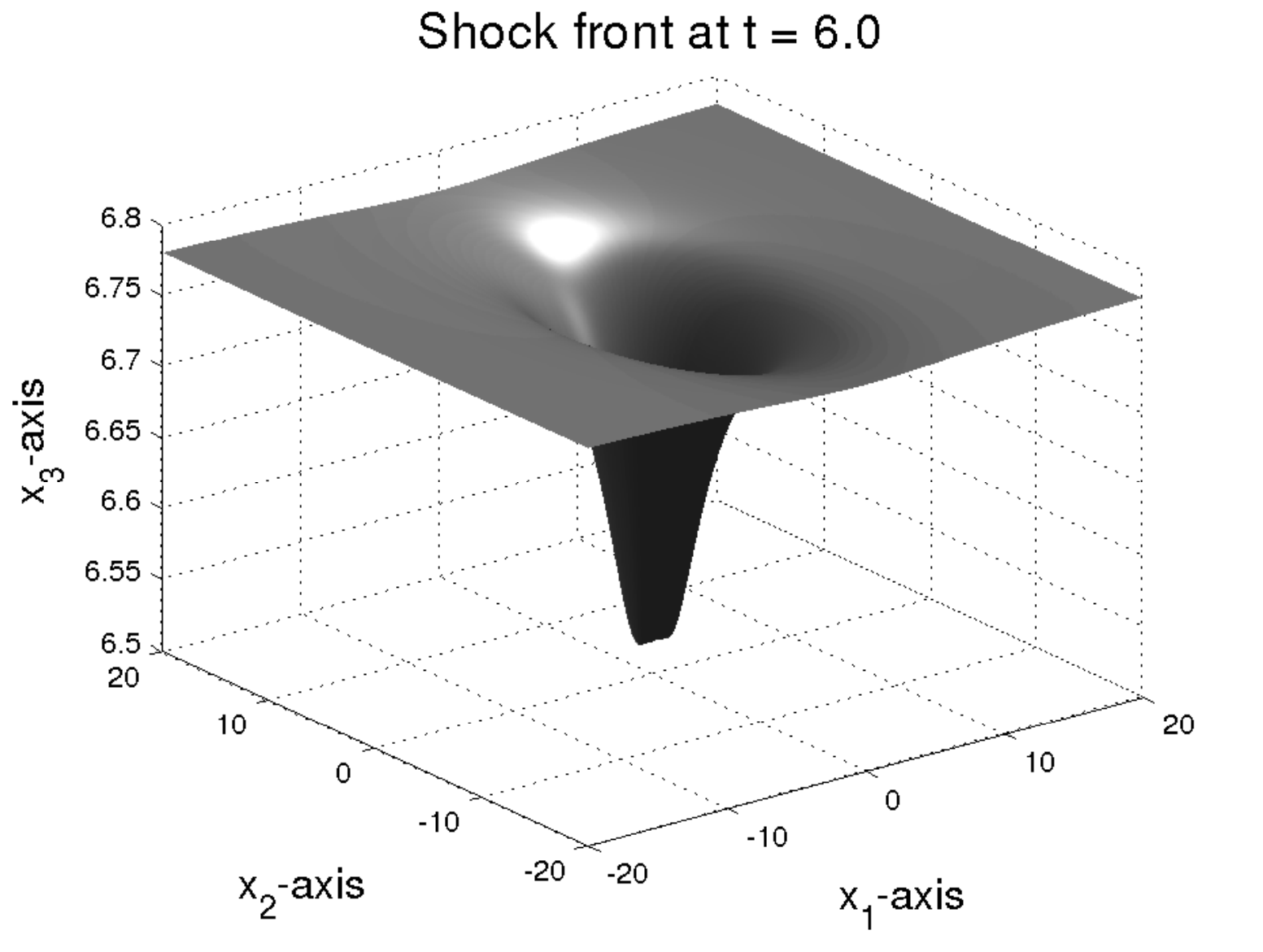}&
    \includegraphics[width=0.5\textwidth]{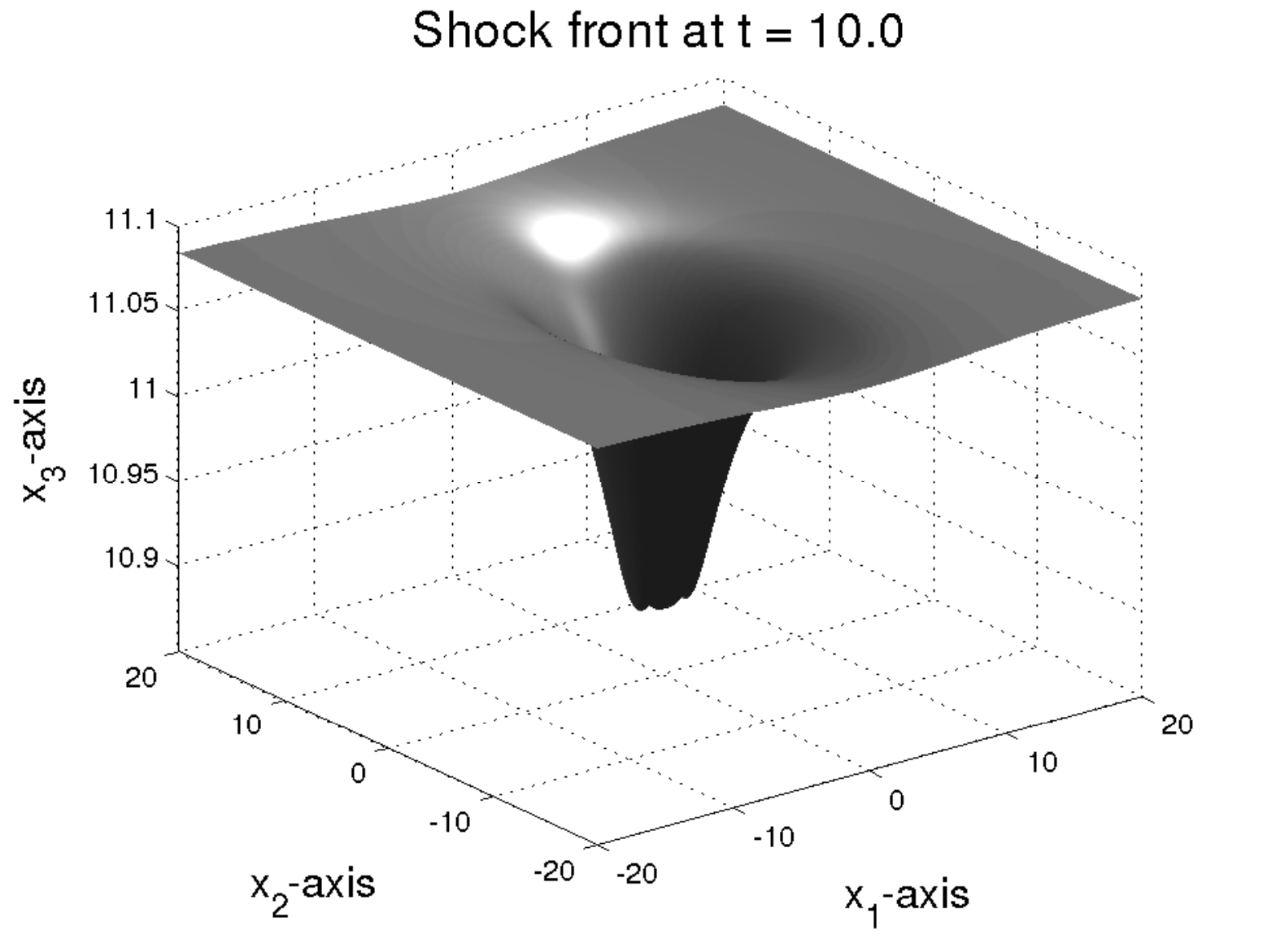}
  \end{tabular}
  \caption{The successive positions of the shock front $\Omega_t$ with
    an initial smooth non-symmetric dip which is not axisymmetric.}
  \label{dip_srt_3d}
\end{figure}

\begin{figure}[htbp]
  \centering
  \includegraphics[width=0.55\textwidth]{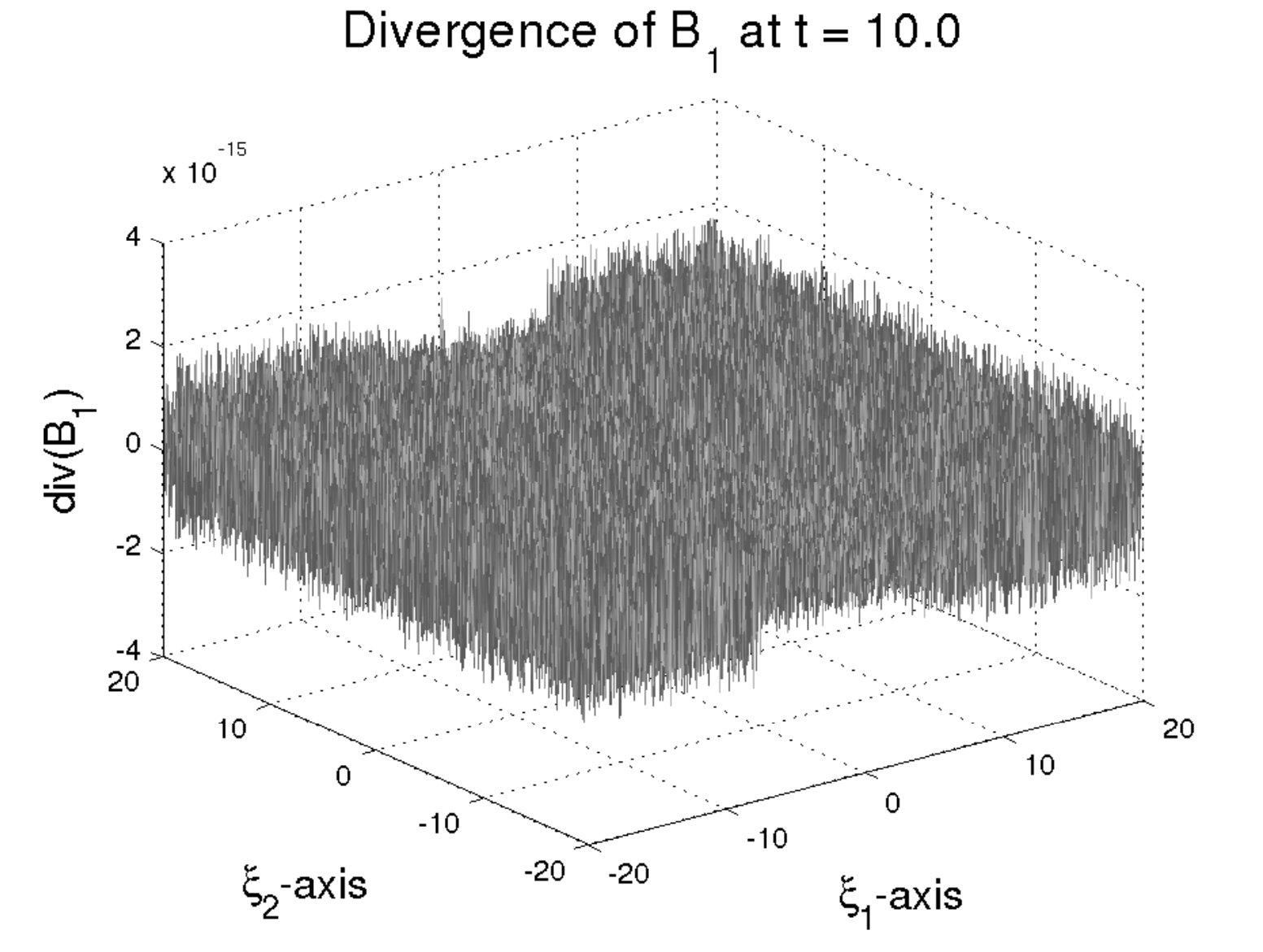}
  \caption{The discrete divergence of $\mathfrak{B}_1$ at
    $t=10.0$. The error is of the order of $10^{-15}$.}
  \label{dip_srt_div}
\end{figure}

In Figure~\ref{dip_srt_3d} we plot the initial shock front
$\Omega_0$ and the successive positions of the shock front
$\Omega_t$ at times $t=2.0,6.0,10.0$. It can be seen that the
whole shock front has moved up in the $x_3$-direction and the dip
has spread over a larger area in $x_1$- and $x_2$-directions. Due
to focusing, the velocity at the central convex portion increases
and as a result, this part of the front moves up, leading to a
change in shape of the initial front. However, when comparing with
the results of 3-D WNLRT in \cite{arun-ct}, we notice that
vertical movement of the corresponding nonlinear wavefront is more
and its central portion bulges out, which has not happened for the
shock front till $t=10.0$. This different behaviour of a shock
front is due to the decay of the shock amplitude as a result of
the nonlinear waves catching-up from behind, represented
mathematically by the variable $\mathcal{V}$.


\subsection{Corrugation Stability of a Shock Front and Interaction
  of Kink Lines}
\label{subsec:corrug}

\begin{figure}
  \centering
  \includegraphics[height=0.6\textheight]{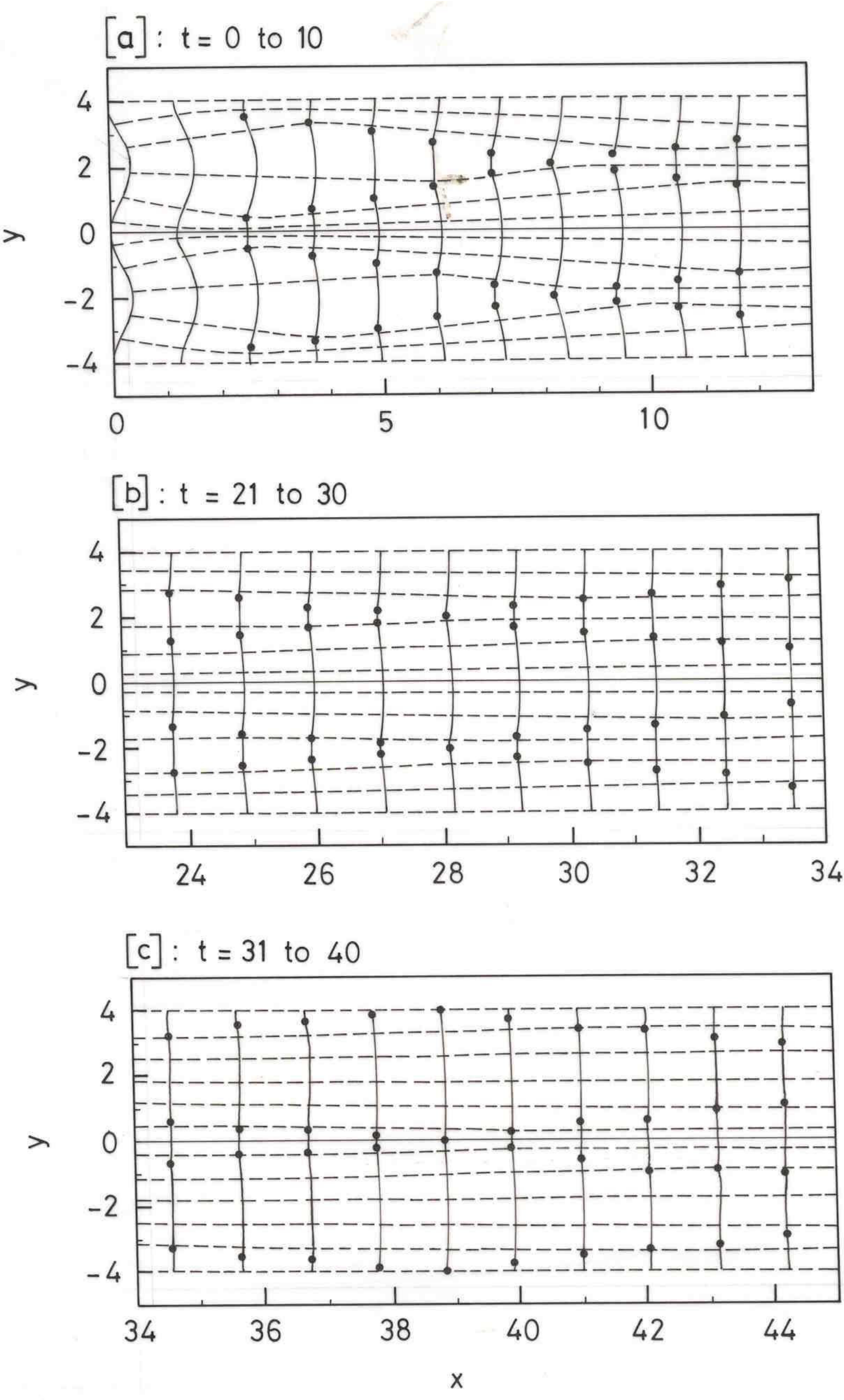}
  \caption{Successive positions of an initially sinusoidal shock front
    and rays from \cite{monica}, plotted at $t=0,1,2,3,\ldots,40$. The
    initial shock front is $x=0.2-0.2\cos(\pi y/2)$ with $M_0=1.2$ and
    $\mathcal{V}_0=0.1$.}
  \label{2d_srt_sine}
\end{figure}

The corrugation stability of a front is defined as the stability
of a planar front to perturbations. This means that the
perturbations in the shape of a planar front ultimately disappear
as time tends to infinity. The extensive numerical simulations by
Monica and Prasad \cite{monica}, using 2-D SRT, clearly show that
a 2-D shock front is corrugation stable.  From \cite{monica} we
reproduce Figure~\ref{2d_srt_sine}, where continuous lines are
successive positions of the 2-D initially sinusoidal shock front,
the broken lines are the rays and dots are the kinks. The shock
has become almost a straight line much before $t=40$. Similarly,
the results of numerical experiments with 3-D WNLRT, reported in
\cite{arun-ct,arun-etal-siam}, show that a 3-D nonlinear wavefront
is also corrugation stable. The aim of this test case is to verify
the corrugation stability of a 3-D shock front, evolving according
to 3-D SRT, and describe the interaction of kink
lines.

The corrugation stability is a result of the genuine nonlinearity
in the characteristic fields corresponding to the two nonzero
eigenvalues of the system (\ref{3d_srt_cons}). The shocks in the
$(\xi_1,\xi_2,t)$-coordinates, which are mapped onto kinks, cause
dissipation of the kinetic energy. Notice that the energy
transport equation \eqref{3d_wnlrt_cons3} of 3-D WNLRT is
homogeneous, whereas the corresponding equation
\eqref{3d_srt_cons3} of 3-D SRT has a source term. In the case of
a nonlinear wavefront, the value of $m-1$ converges to the mean
value of $m_0-1$. For a shock, when $\mathcal{V}>0$, the value of
$M-1$ decreases to zero. This is typical of a plane shock in gas
dynamics, which can also be seen from the 1-D model equation
$u_t+(u^2/2)_x=0$; see \cite{prasad-book}. Thus, the value of
$M-1$, in addition to approaching a constant value, decays. The
result is that the perturbations in the shape of the shock front
not only disappear, leading to corrugation stability, but also the
front velocity $M$ approaches the linear front velocity $M=1$. We
refer the reader to \cite{monica} for a related discussion on the
corrugation stability of 2-D shock fronts.

\begin{figure}
  \centering
  \includegraphics[width=0.5\textwidth]{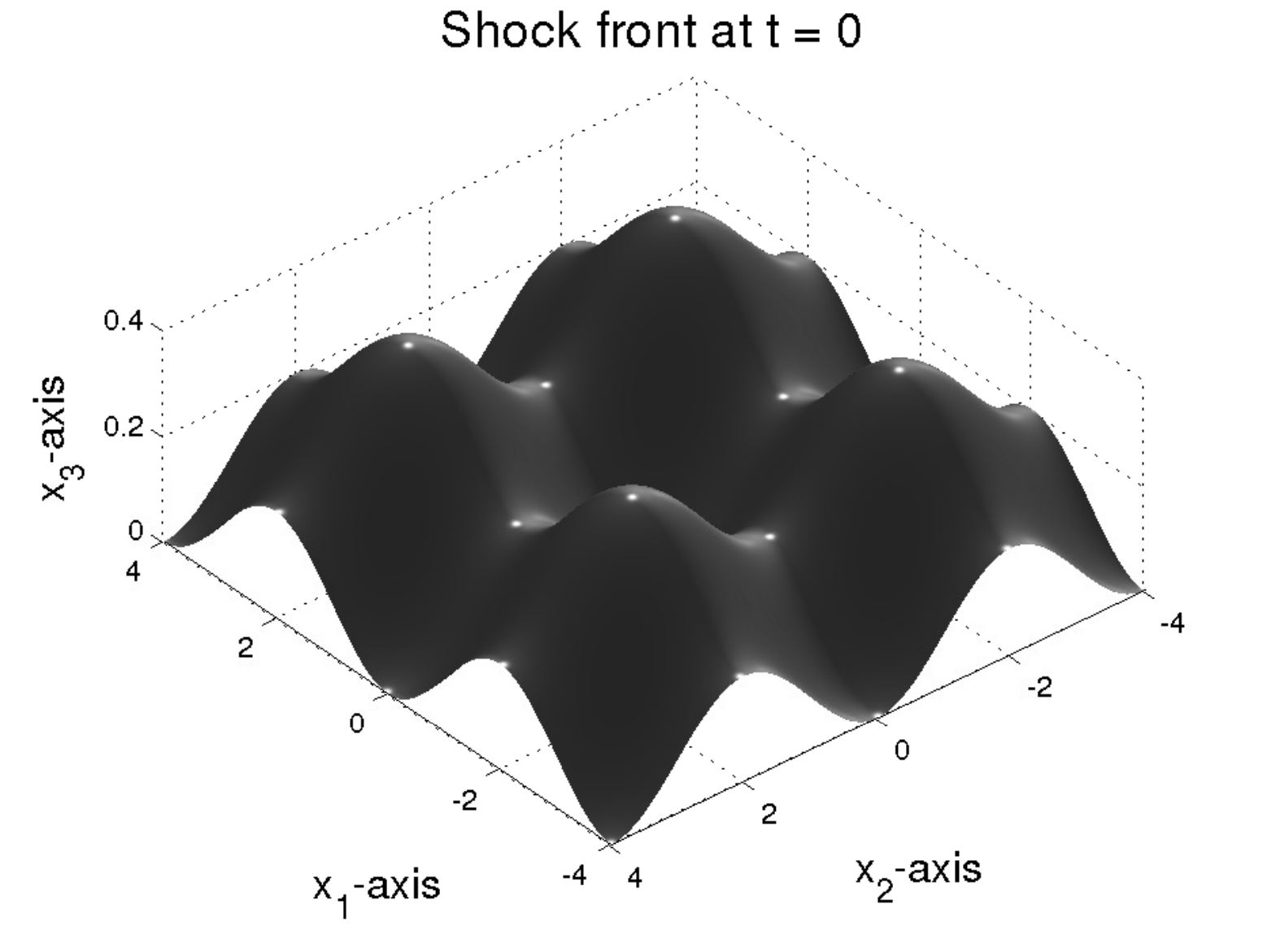}
  \caption{Initial shock front in the shape of a smooth periodic
    pulse.}
  \label{sine_srt_3d_0}
\end{figure}

In order to verify the corrugation stability, we consider here the
initial shock front $\Omega_0$ to be of a periodic shape in $x_1$- and
$x_2$-directions
\begin{equation}
  \Omega_0\colon x_3=\kappa\left(2
    -\cos\left(\frac{\pi x_1}{a}\right)-\cos\left(\frac{\pi
        x_2}{b}\right)\right)
  \label{sine_omega0}
\end{equation}
with the constants $\kappa=0.1,a=b=2$. In
Figure~\ref{sine_srt_3d_0} we give the plot of the initial shock
front $\Omega_0$, which is a smooth pulse without any kink lines.
The initial front $\Omega_0$ can be thought of modelling a smooth
perturbation of a plane front. We prescribe a constant initial
velocity $M_0=1.2$ everywhere on the front given in
(\ref{sine_omega0}). Though the initial shock front is smooth, as the
time evolves, a number of kink lines appear in each period. The
process of interaction of these kink lines is a very interesting
phenomena, which has not been described in \cite{arun-ct}. In the
following, we proceed to do it with a number of plots of the shock
front at different instances.

\begin{figure}
  \centering
  \begin{tabular}{cc}
    \includegraphics[width=0.5\textwidth]{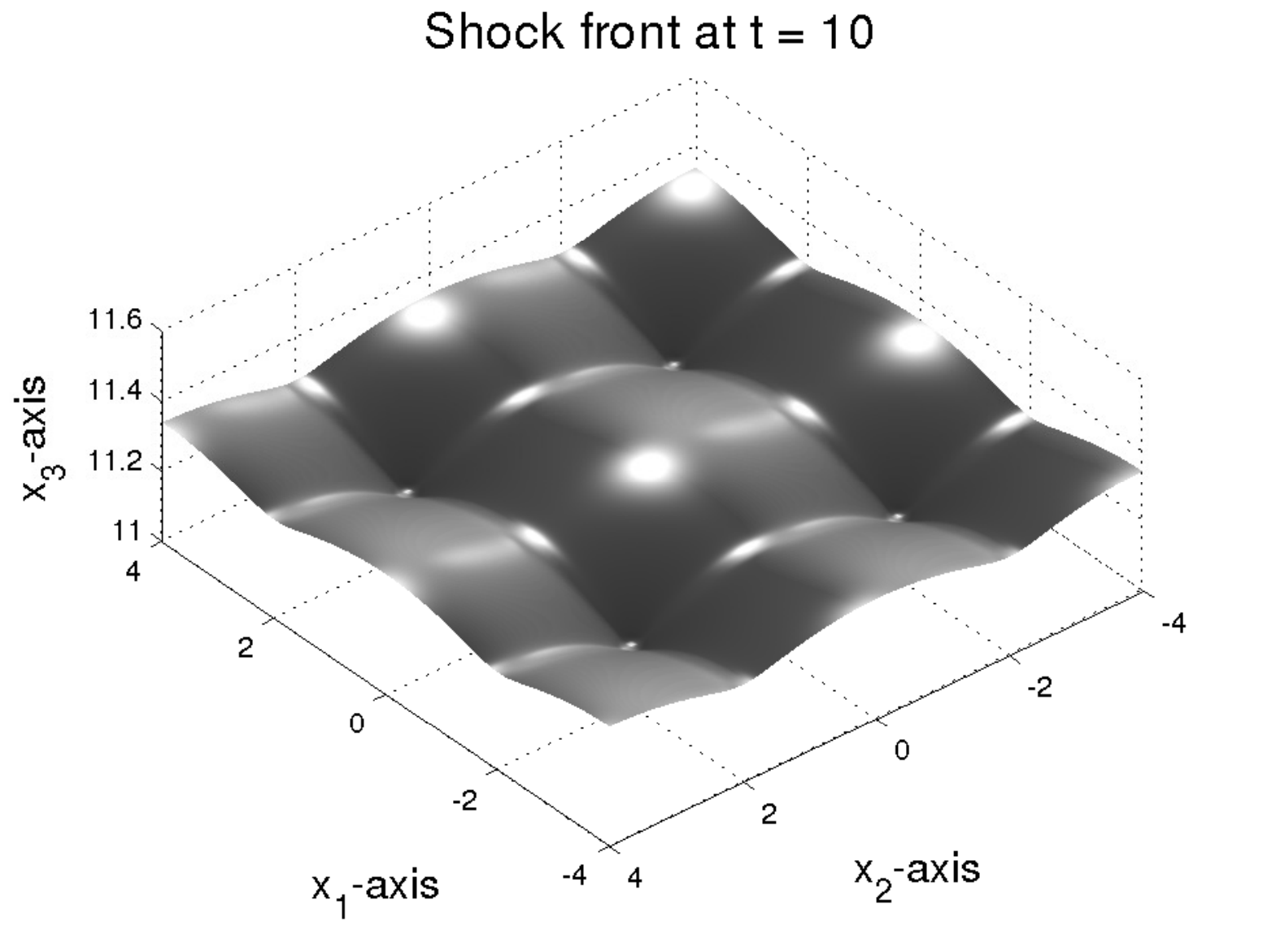}&
    \includegraphics[width=0.5\textwidth]{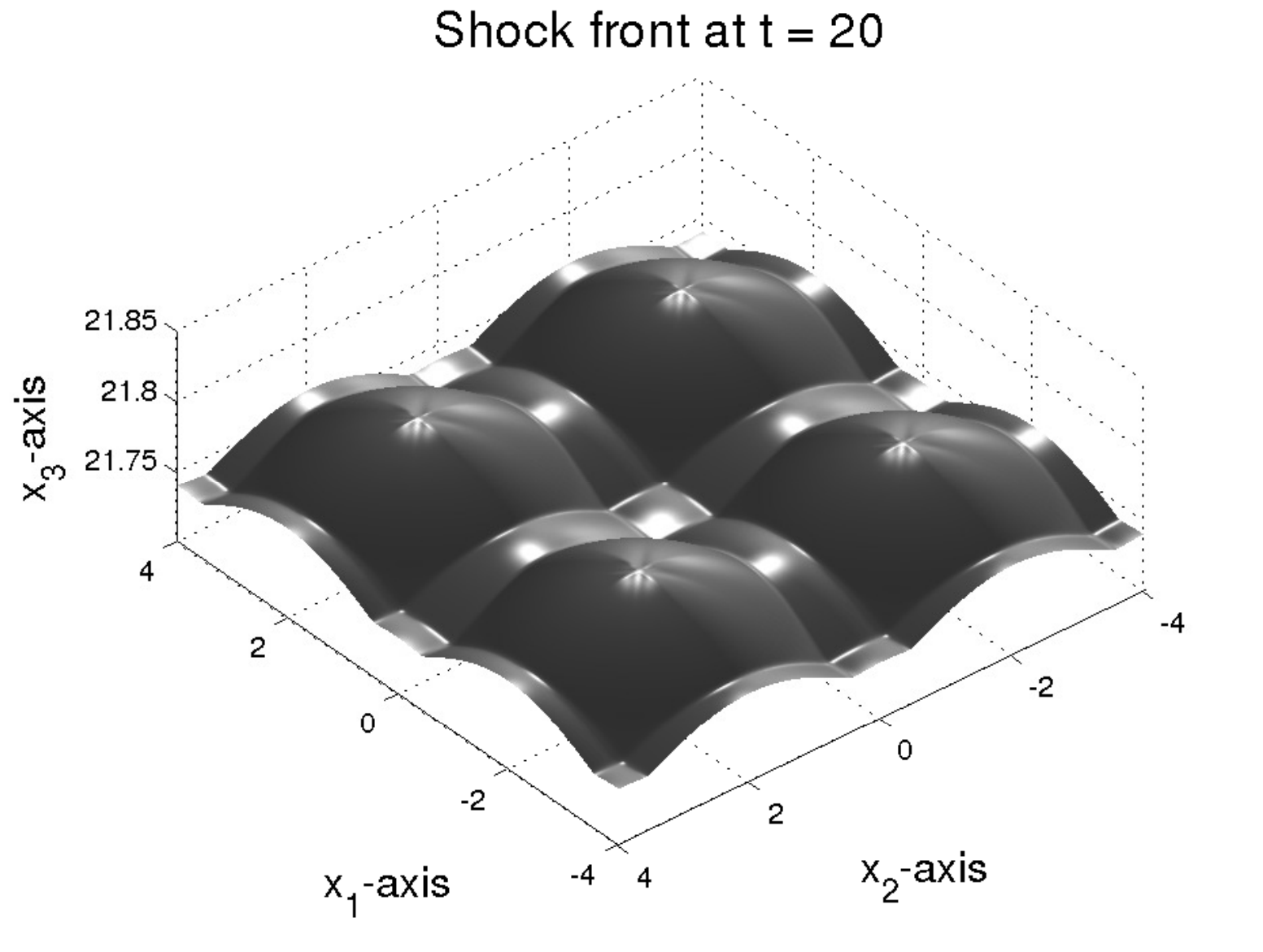}\\
    \includegraphics[width=0.5\textwidth]{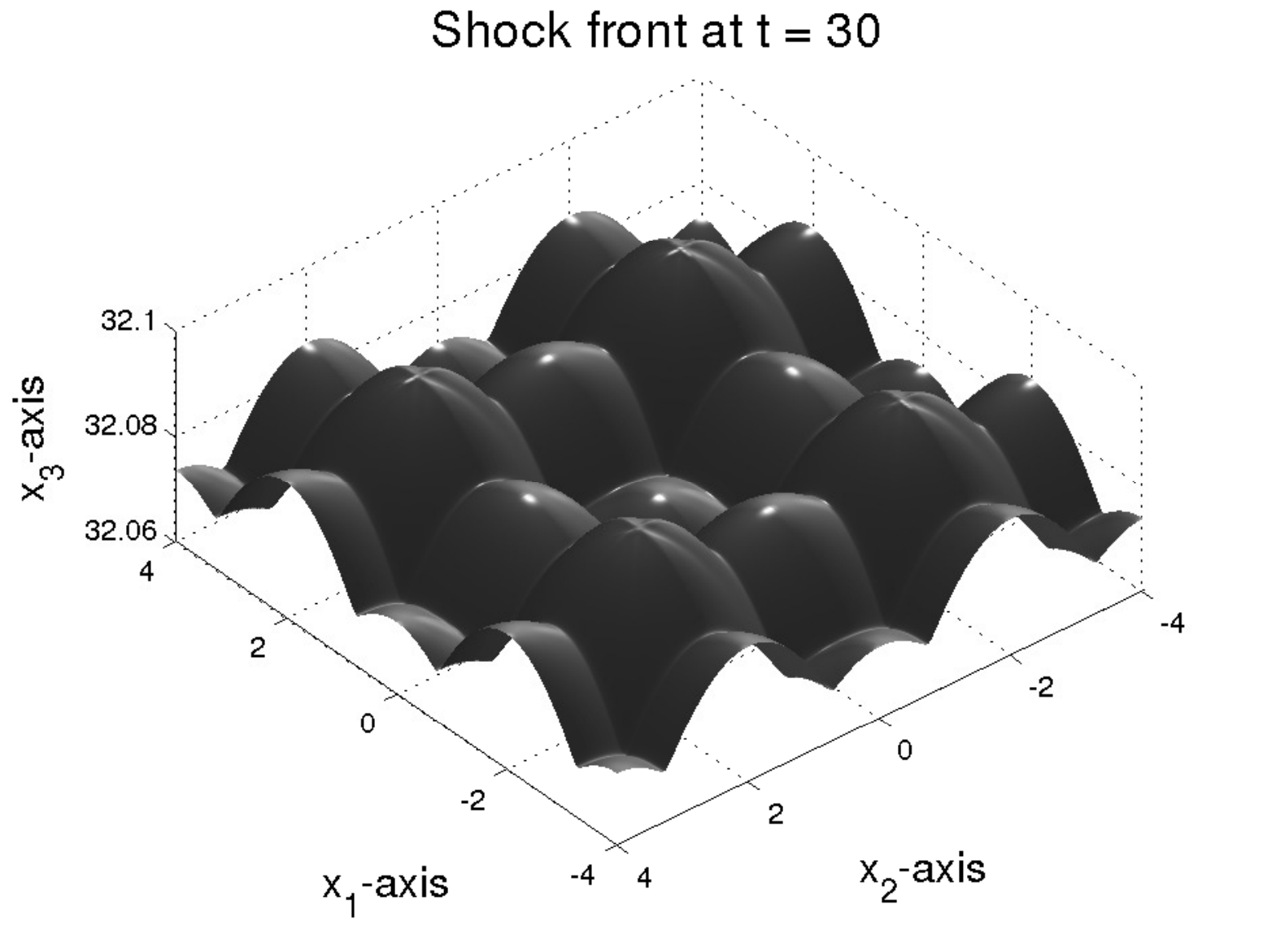}&
    \includegraphics[width=0.5\textwidth]{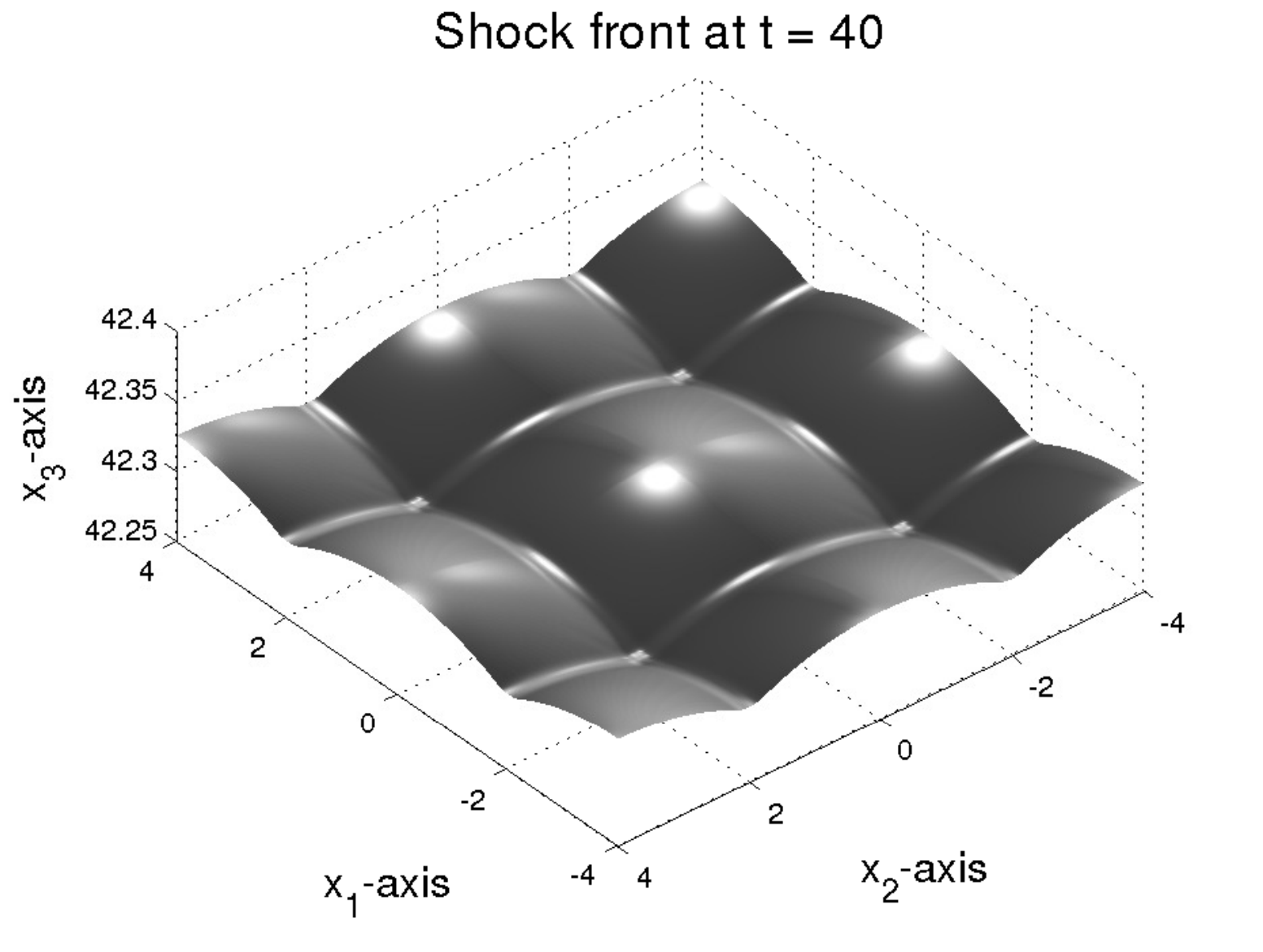}\\
    \includegraphics[width=0.5\textwidth]{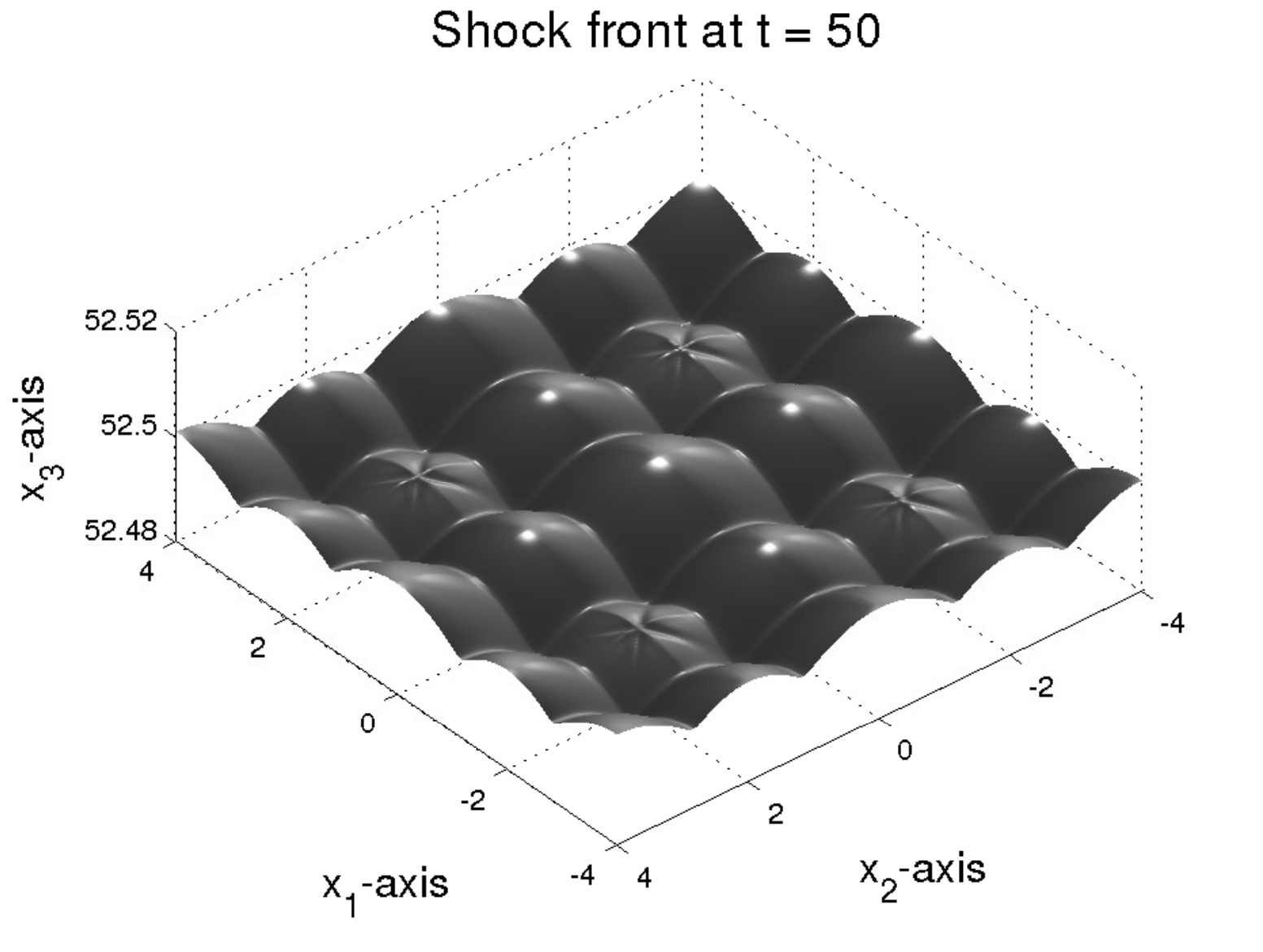}&
    \includegraphics[width=0.5\textwidth]{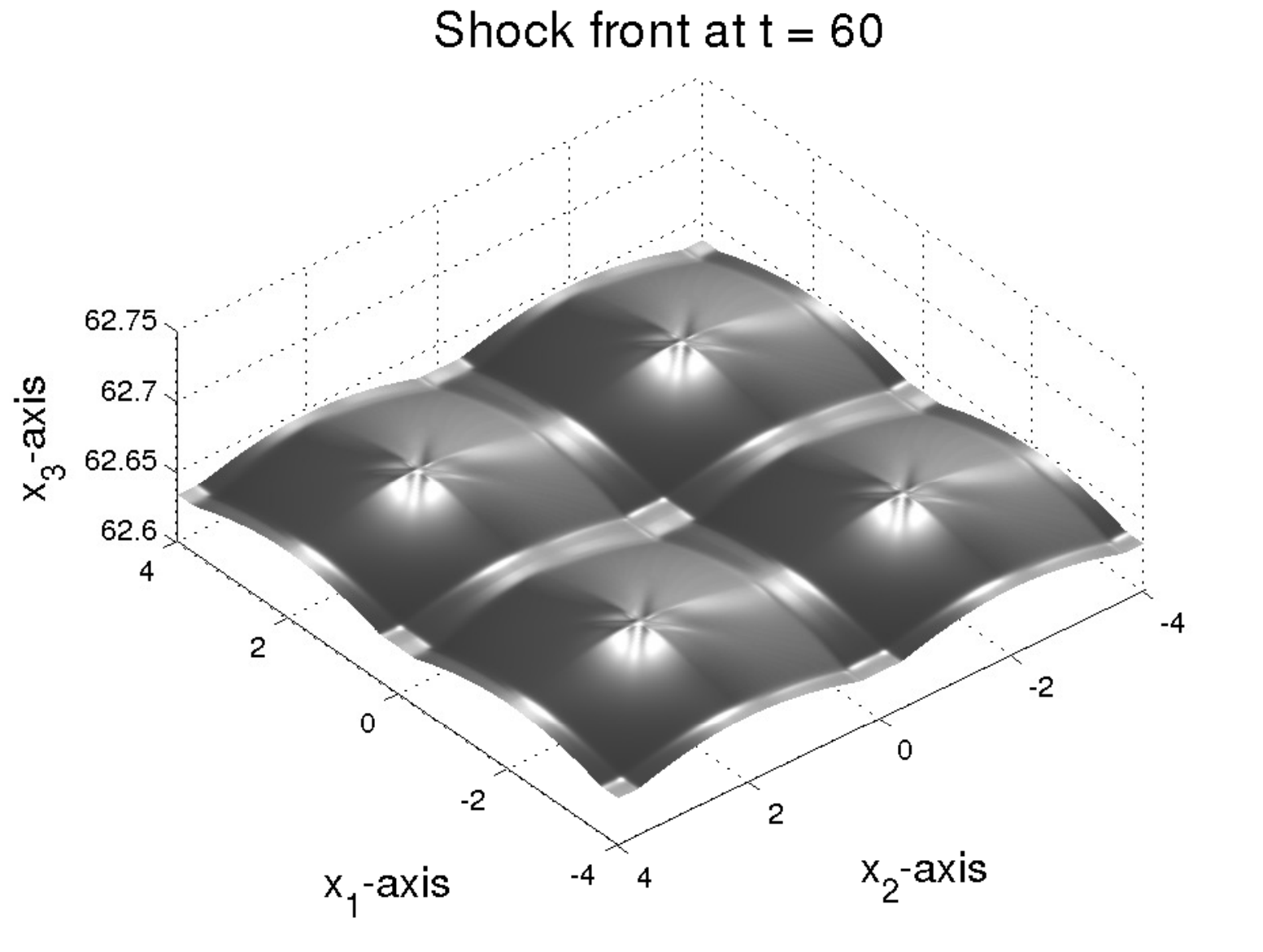}
  \end{tabular}
  \caption{Shock front $\Omega_t$ starting initially in a periodic
    shape with $M_0=1.2$. The shock front develops a complex pattern
    of kinks and ultimately becomes planar.}
  \label{sine_srt_3d}
\end{figure}

\begin{figure}[htbp]
  \centering
  \includegraphics[width=0.45\textwidth]{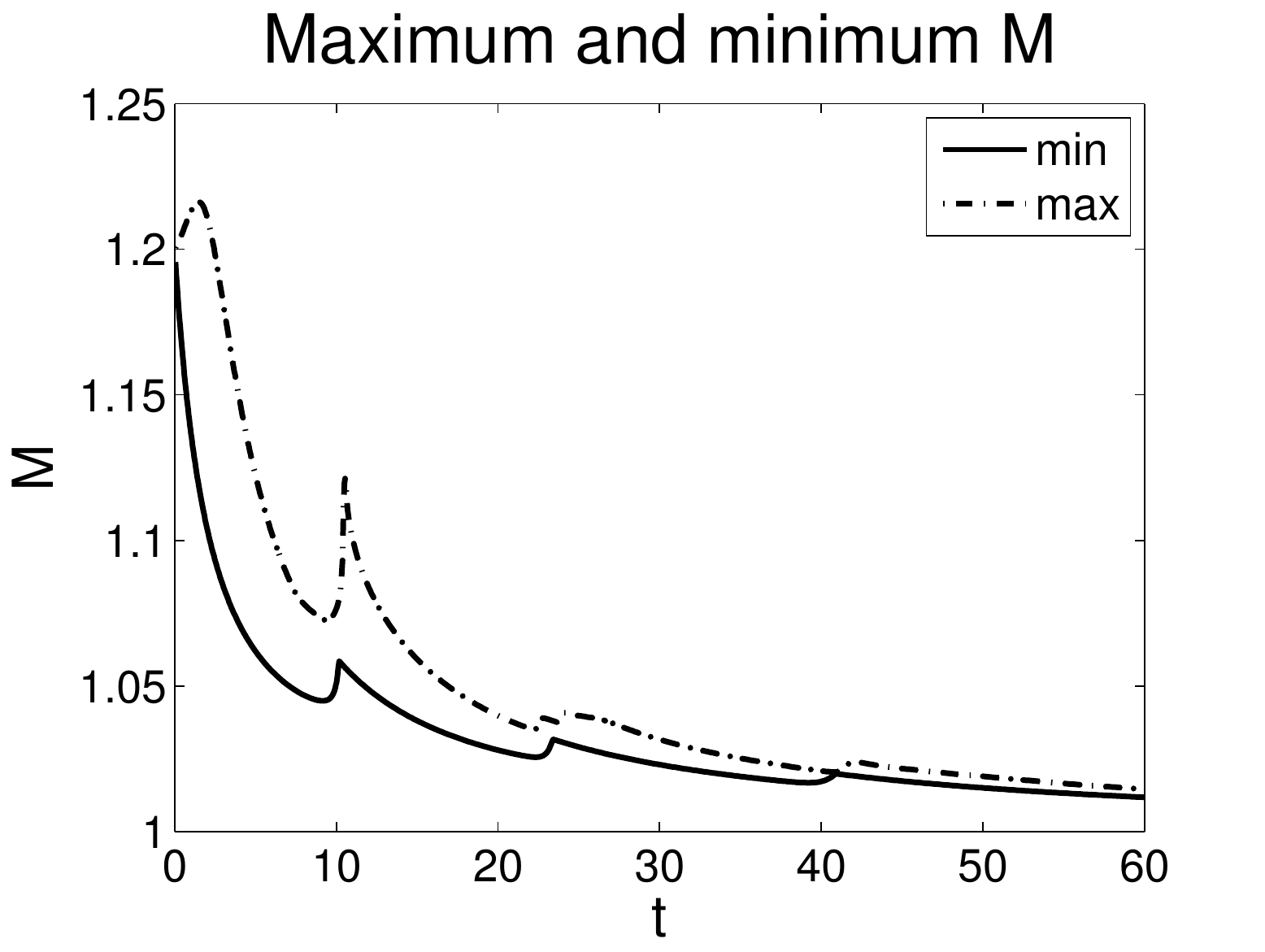}
  \caption{Variation of $M_{\max}(t)$ and $M_{\min}(t)$ with time
    from $t=0$ to $t=80$ for a periodic shock front $x_3=\kappa\left(2
    -\cos\left(\frac{\pi x_1}{a}\right)-\cos\left(\frac{\pi
        x_2}{b}\right)\right)$. The difference $M_{\max}(t)-M_{\min}(t)$
  tends to zero as $t\to\infty$.}
  \label{sine_srt_min_max}
\end{figure}

In Figure~\ref{sine_srt_3d} we give the surface plots of the shock
front $\Omega_t$ at times $t=10,20,30,40,50,60$ in two periods in
$x_1$ and $x_2$-directions. As mentioned above, the initial shock
front is smooth, with no kink lines. The front $\Omega_t$ moves up in
the $x_3$-direction and develops several kink lines. Four kink lines
parallel to $x_1$-axis and four parallel to $x_2$-axis can be seen in
the figures on the shock front at times $t\geq 10$. These kink lines
are formed at a time before $t=10$, say about $t=2$. This can easily
be observed from the Figure~\ref{sine_srt_min_max} showing the maximum
and minimum values, $M_{\max}(t)$ and $M_{\min}(t)$, versus $t$,
where $M_{\max}(t)$ and $M_{\min}(t)$ the maximum and minimum values
of $M$ respectively taken over $(\xi_1,\xi_2)$ at any time $t$. We can
observe a significant increase in the maximum value of $M$ near $t=2$
and it also jumps after the interaction of kink lines at later times.


Due to symmetry, it is sufficient if we describe the motion of the
kink lines parallel to the $x_1$-axis. Let us designate the part
$-4\leq x_1\leq 0$ as the first period and $0 \leq x_1 \leq 4$ as
second period. A pair of kink lines are seen in each period at $t=10$
near $x_1=-2$ and $x_1=2$, they are about to interact and produce
another pair of kink lines. These newly formed kink lines move apart,
cf.\ also the corresponding 2-D diagram in
Figure~\ref{2d_srt_sine}. At $t=20$, one kink line from each of first
and second periods have come quite near and are seen close to
$x_1=0$. They interact, produce a new pair of kink lines and move
apart as seen at $t=30$, where we find four distinct kink lines in
$-4<x_1<4$. This process continues, two pairs of kink lines about to
interact are seen near $x_1=-2$ and $x_1=2$ at $t=40$. There are four
distinct kink lines also at $t=50$ and at $t=60$. The two kink lines
near $x_1=0$  at $t=60$ are about to interact.

The interaction of kinks was first observed on a 2-D shock front
experimentally in \cite{sturtevant-kulkarny}, cf.\ Figures 6(b) and
(c) and numerically in \cite{monica}, cf.\
Figure~\ref{2d_srt_sine}. We do observe the same phenomenon here,
except that kinks are replaced by kink lines. A detailed theoretical
analysis of this phenomenon on a 2-D nonlinear wavefront is presented
in \cite{baskar-riemann}. Interaction of kinks is a particular case of
some interesting phenomena and in the following we briefly describe
them; see also \cite{baskar-riemann} for more details.

The equations of the 2-D KCL based WNLRT form a hyperbolic system of
three conservation laws in $(\xi_1,t)$-plane when $m>1$; see also
\eqref{3d_wnlrt_cons1}-\eqref{3d_wnlrt_cons3} and \eqref{eq:m_mu}.
The system of conservation laws of 2-D WNLRT is not degenerate and
there are no source terms. The characteristic fields corresponding to
the two nonzero eigenvalues are genuinely nonlinear and corresponding
to these fields we have four elementary wave solutions, which are
shocks and centred rarefaction waves in $(\xi_1,t)$-plane. The images
of these elementary waves onto $(x_1,x_2)$-plane are called
elementary shapes. As mentioned earlier, a shock is mapped onto an
elementary shape kink in $(x_1,x_2)$-plane but a centred rarefaction
wave is mapped on a continuous convex shape on the weakly nonlinear
wavefront $\Omega_t$. Interaction of the elementary shapes in
$(x_1,x_2)$-plane can be studied with the help of interactions of
elementary waves in $(\xi_1,t)$-plane. These interactions can be of
finite or infinite duration in time. Two kinks of different
characteristic family on $\Omega_t$ approach, interact and then
produce another pair of kinks, which move apart.

For $M>1$, the KCL based 2-D SRT equations form a hyperbolic
system of balance laws due to appearance of second terms in
\eqref{2d_srt_cons3}-\eqref{2d_srt_cons4} which are source
terms. However, the source terms, which affect the solution on larger
space and time scales, do not have any effect on the interaction of
two kinks on a shock front, which takes place instantaneously. They
have only a small effect on the motion of kinks for a short time
before and after the interaction. Thus, the nature of interaction of
two kinks on a shock front as seen in Figure~\ref{2d_srt_sine} is same
as that of interaction of kinks on a nonlinear wavefront theoretically
predicted in \cite{baskar-riemann}. Interactions of two parallel
kink lines on a 3-D shock front will be similar to that of the two
kinks on a 2-D shock front. We clearly observe this in
Figure~\ref{sine_srt_3d} for the interaction of a pair of parallel
kink lines. However, the interaction of oblique kink lines will be
quite different about which we do not intent to make any comment.

\begin{remark}
  Comparing the shock fronts at different times in
  Figure~\ref{sine_srt_3d}, with the corresponding nonlinear
  wavefronts at the same time, we notice that the two fronts differ
  in their shapes. Let us compare the graphs of $M_{\max}(t)$ and
  $M_{\min}(t)$ in Figure~\ref{sine_srt_min_max} with those of
  $m_{\max}(t)$ and $m_{\min}(t)$ in Fig.11(a) of \cite{arun-ct}. We find
  that $M_{\max}(t)$ and $M_{\min}(t)$ values on the shock front decay
  very rapidly compared to those on the nonlinear wavefront. As a
  result, the  waves on the nonlinear wavefront move faster and the
  interaction of kink lines takes place more frequently. Graphs of
  $m_{\max}(t)$ and $m_{\min}(t)$ in Fig.10(a) of \cite{arun-ct}
  oscillate quite fast. Since sudden increases in the values of
  $M_{\max}(t)$ and $M_{\min}(t)$ (and similarly those of
  $m_{\max}(t)$ and $m_{\min}(t)$) correspond to interactions of kink
  lines, we notice that kink lines on the nonlinear wavefront interact
  more frequently than those on a shock front.
\end{remark}

\subsection{A Shock Front Starting from an Axisymmetric Shape}
\label{subsec:axisym}

Next we consider a shock front with an axisymmetric, oscillatory and
radially decaying shape given by
\begin{equation}
  \Omega_0\colon x_3=\kappa\cos(\alpha r)e^{-\beta r},
  \label{cos_exp_omega0}
\end{equation}
where $r=\sqrt{x_1^2+x_2^2}$. The evolution of a nonlinear
wavefront with this initial geometry has not been discussed in
\cite{arun-ct} and we present it here as yet another instance of
corrugation stability. The initial shock front $\Omega_0$ models a
smooth perturbation of a planar front such that the amplitude of the
perturbation decays to zero as $r\to\infty$. The parameters in
\eqref{cos_exp_omega0} are taken as
$\kappa=0.05,\alpha=1.0,\beta=0.15$. The initial velocity has a
constant value $M_0=1.2$ everywhere on the shock front.

\begin{figure}
  \centering
  \begin{tabular}{cc}
    \includegraphics[width=0.5\textwidth]{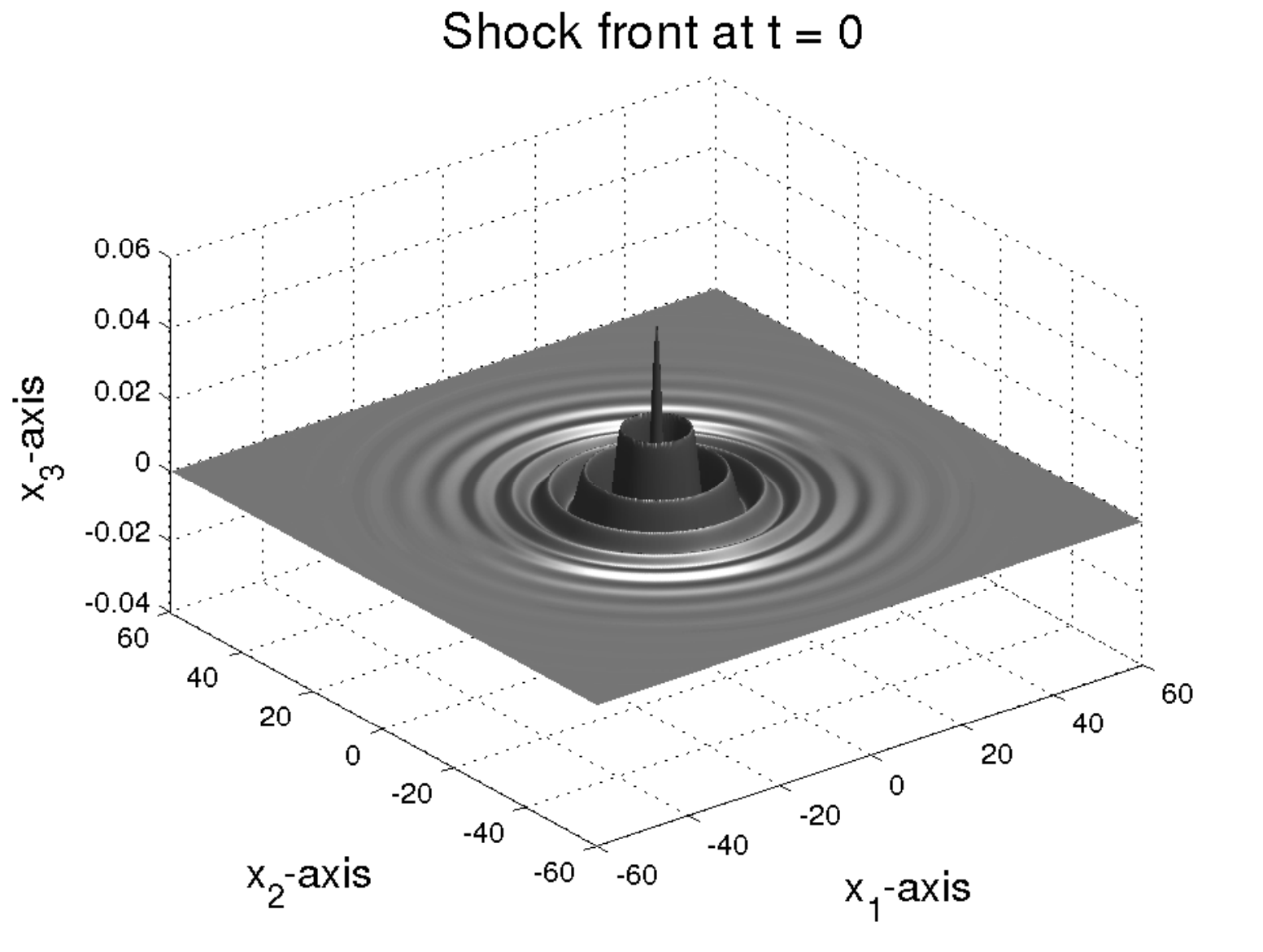}&
    \includegraphics[width=0.5\textwidth]{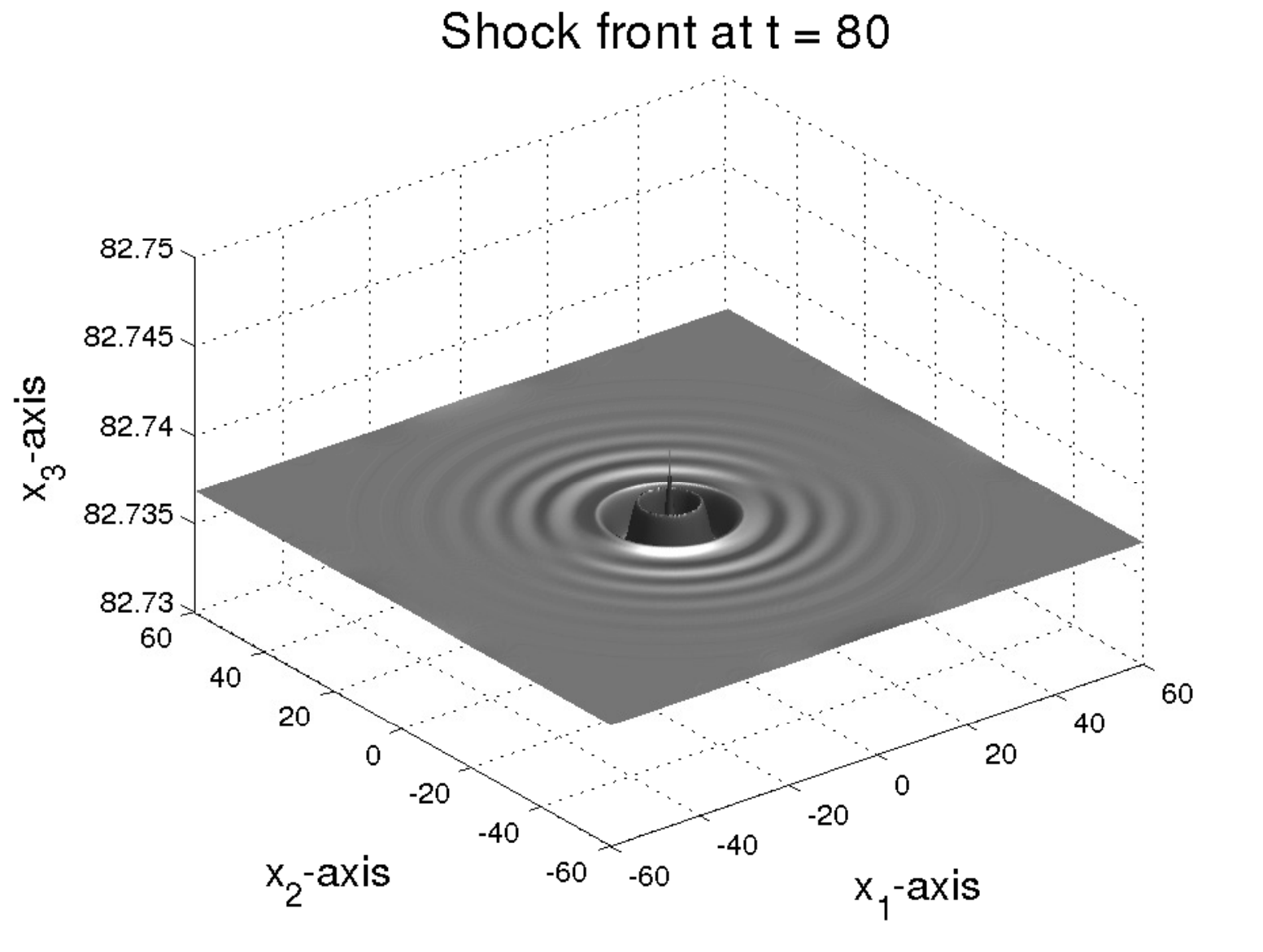}
  \end{tabular}
  \caption{The evolution of a shock front $\Omega_t$ starting from a
    smooth pulse {\bf $x_3=0.05 \cos (r)e^{-0.15r}$}.}
  \label{cos_srt_3d}
\end{figure}

\begin{figure}
  \centering
  \includegraphics[width=0.45\textwidth]{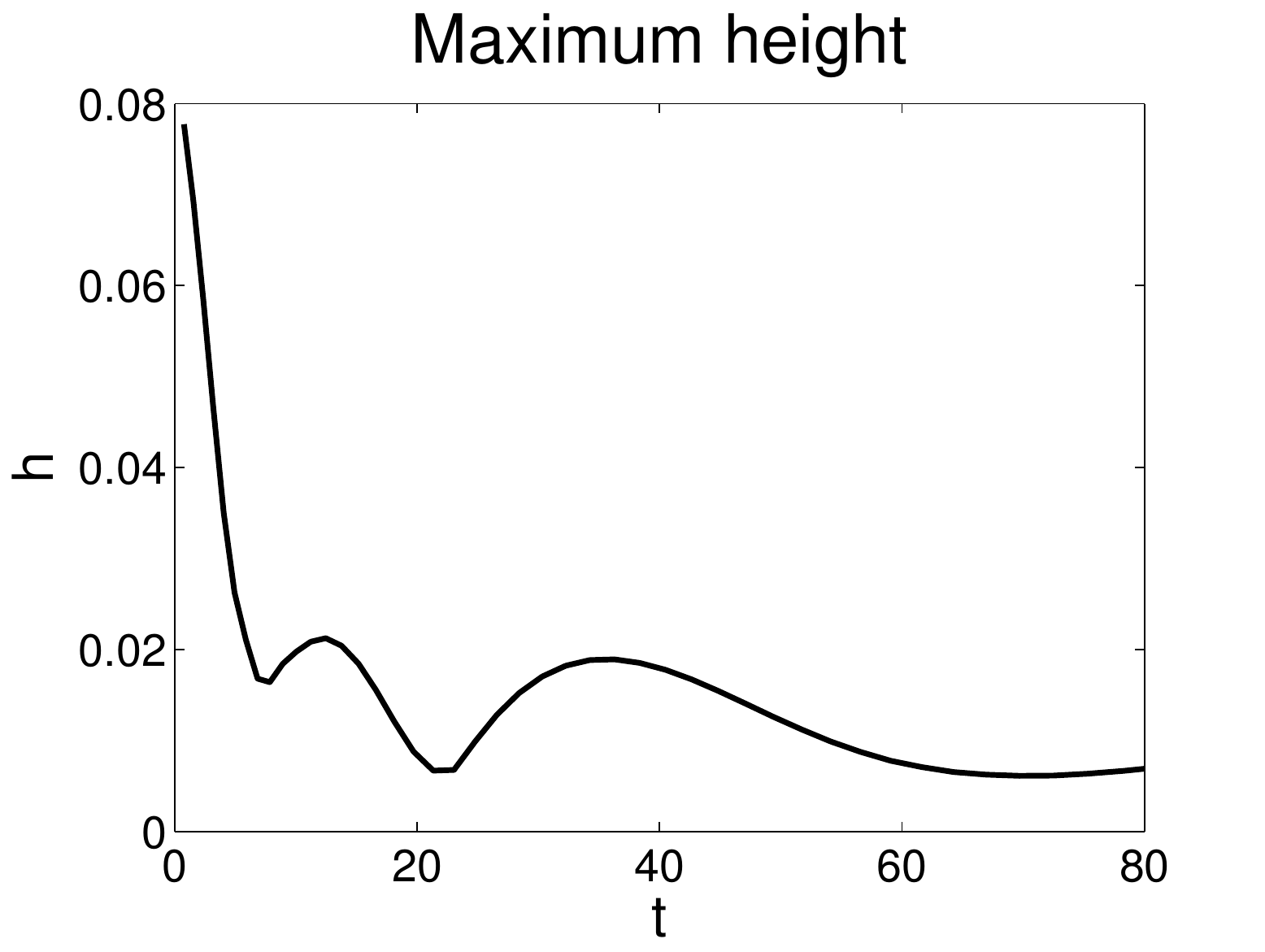}
  \caption{Time variation maximum height of the shock front given
    initially by $x_3=0.05 \cos (r)e^{-0.15r}$.}
  \label{cos_srt_height}
\end{figure}

In Figure~\ref{cos_srt_3d} we give the surface plots of the initial
shock front $\Omega_0$ and the shock front $\Omega_t$ at time
$t=80$. It can be noted that the front $\Omega_t$ moves up in the
$x_3$-direction. At $t=0$, there is an axisymmetric elevation near
the origin $r=0$ and this central elevation decays fast. The smooth
shape $\Omega_0$ at $t=0$ develops later a number of circular kink lines,
however these kink lines have almost disappeared at $t=80$ since the
height of the shock front has become quite small at this time. The
elevations and depressions on the front diminish, leading to the
reduction in height. We compute the maximum height $h(t)$ defined by
\begin{equation}
 h(t):={x_3}_{\max}(t)-{x_3}_{\min}(t),
 \label{max_height}.
\end{equation}
In Figure~\ref{cos_srt_height} we give the plot of $h$ versus $t$,
which clearly shows that height reduces with time. The initial maximum
height is $h(0)=0.08$, whereas at $t=80$ it is $h(80)=0.002917$, which
corresponds to a $96.35\%$ reduction in the initial height. It is,
therefore, very easy to see that the shock front tends to become
planar, with its height decreasing to zero.

\begin{figure}
  \centering
  \begin{tabular}{cc}
    \includegraphics[width=0.45\textwidth]{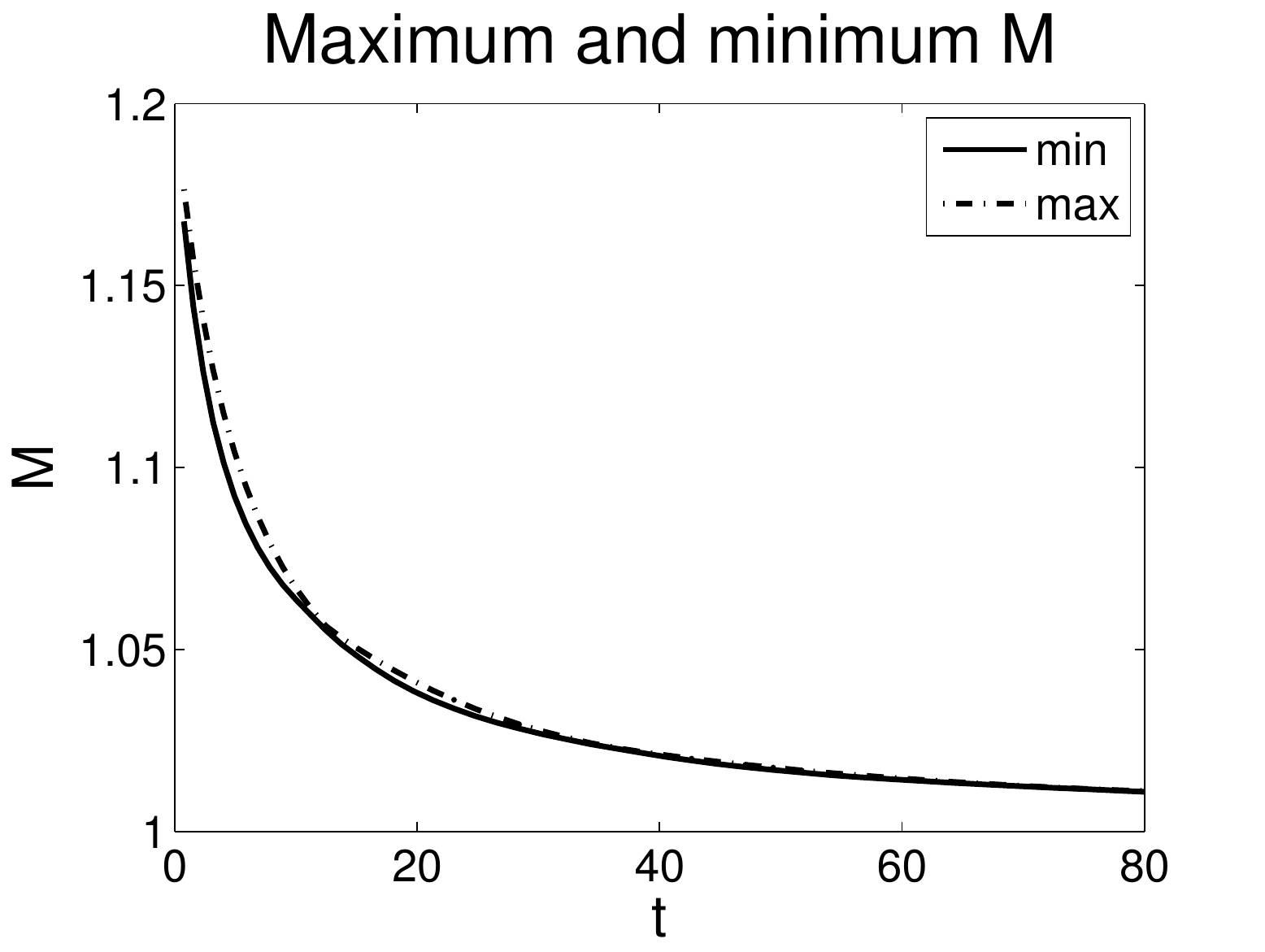}&
    \includegraphics[width=0.45\textwidth]{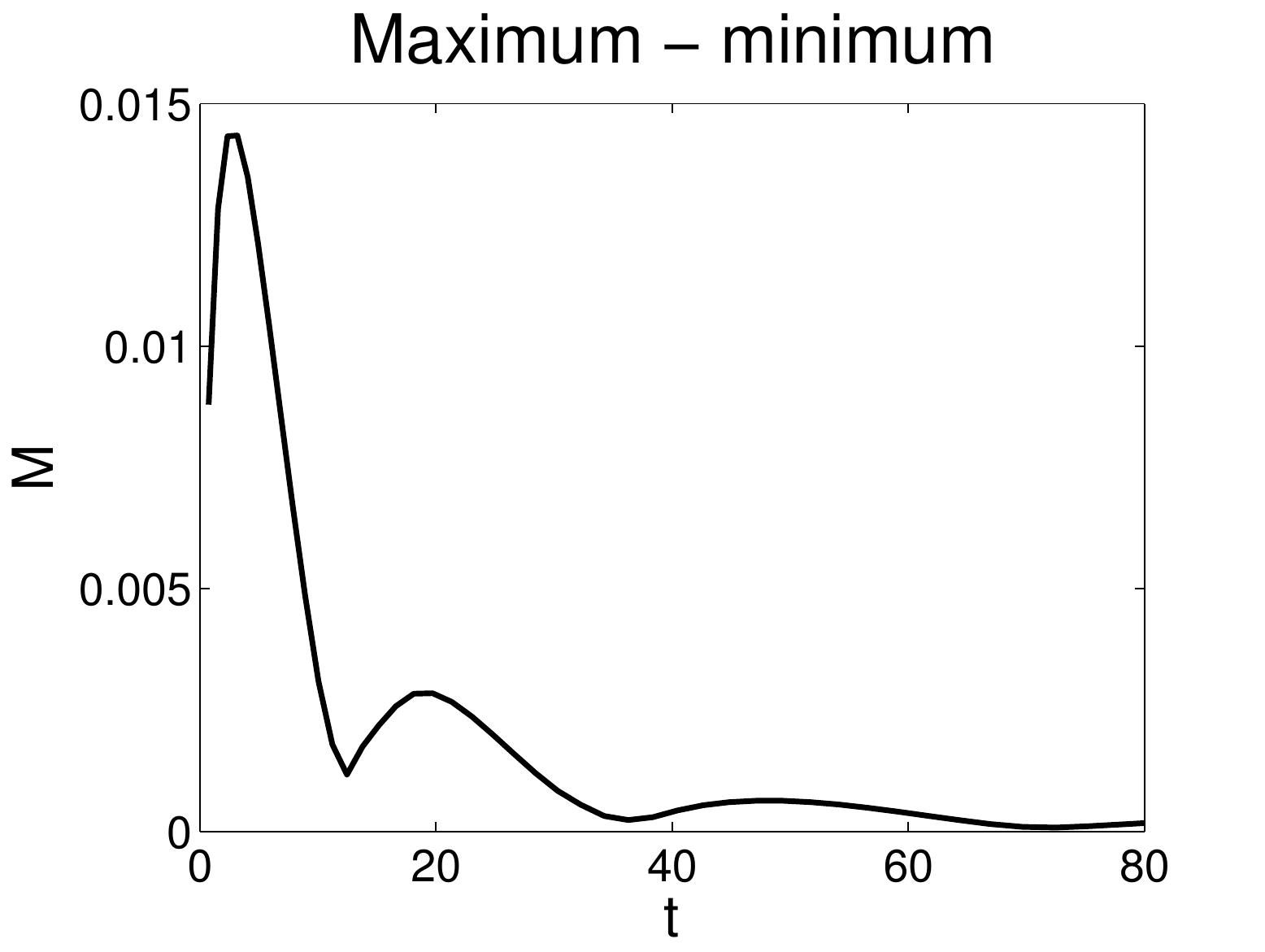}\\
    (a) & (b)
  \end{tabular}
   \caption{For an initial shock front $x_3=0.05 \cos (r)e^{-0.15r}$
     (a): variation of $M_{\max}(t)$ and $M_{\min}(t)$ with time from
     $t=0$ to $t=80$. (b): the difference $M_{\max}(t)-M_{\min}(t)$
     tends to zero as $t\to\infty$.}
   \label{cos_srt_min_max}
\end{figure}

\begin{figure}
  \centering
  \begin{tabular}{cc}
    \includegraphics[width=0.45\textwidth]{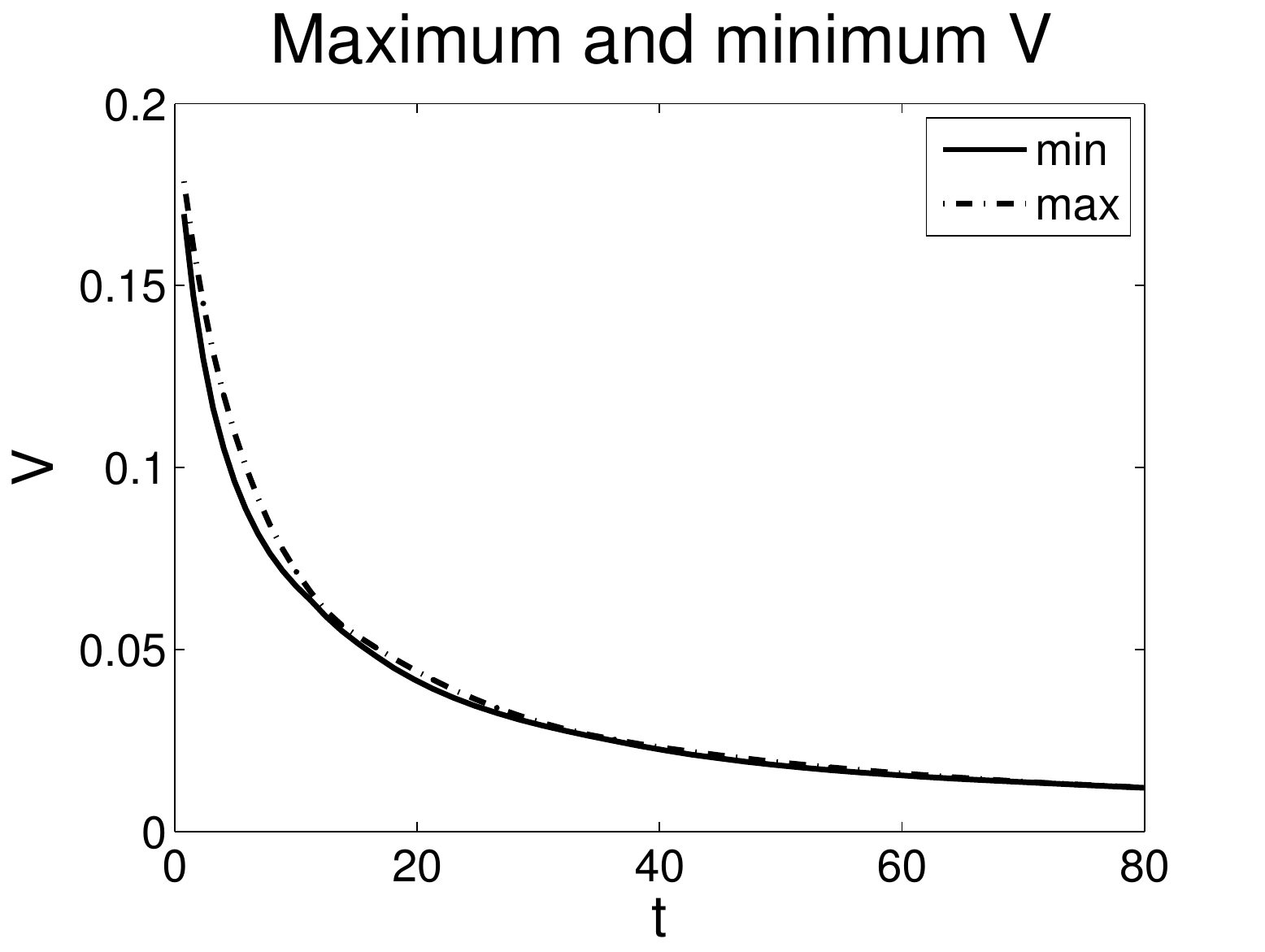}&
    \includegraphics[width=0.45\textwidth]{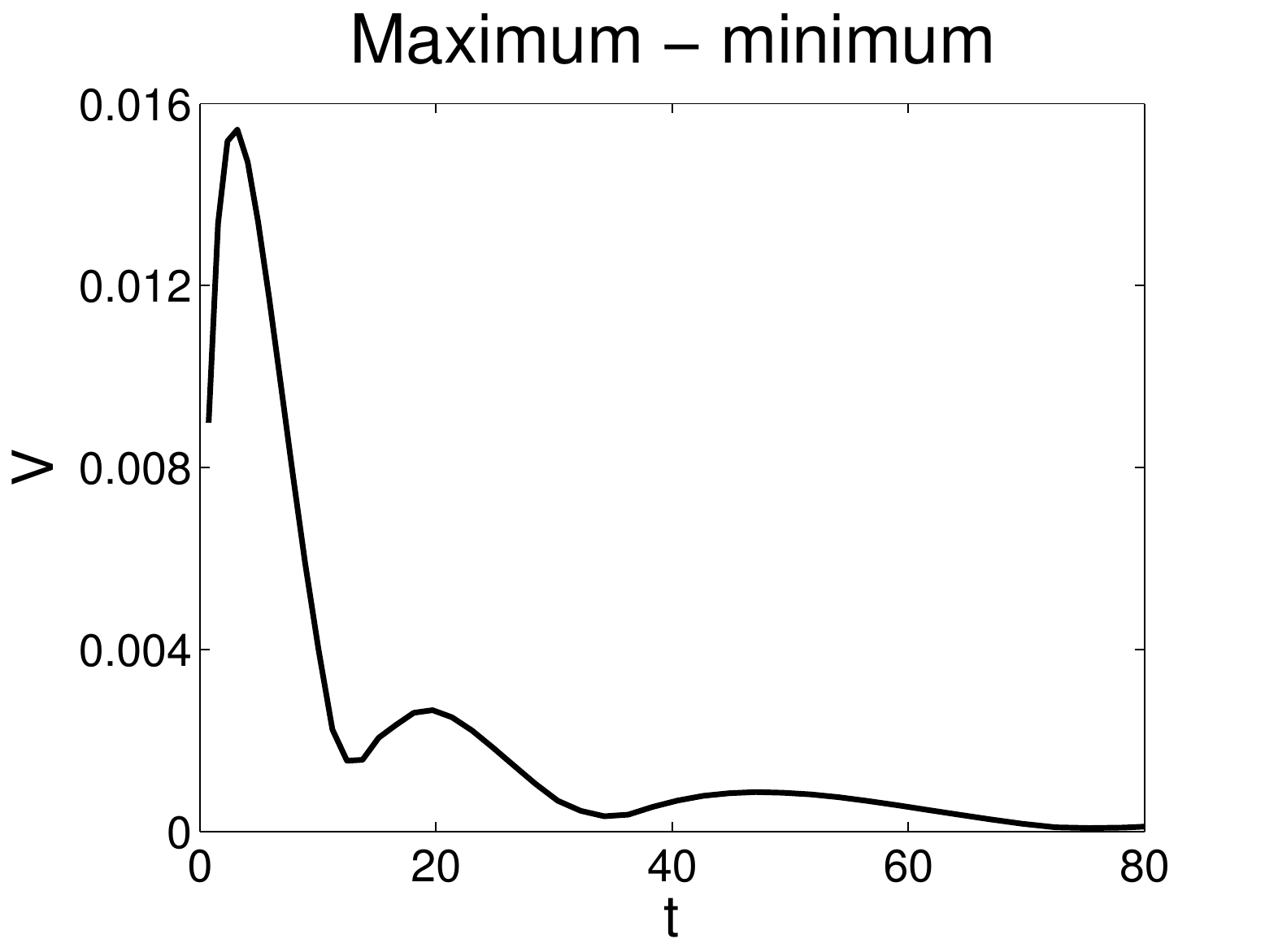}\\
     (a) & (b)
  \end{tabular}
  \caption{For an initial shock front $x_3=0.05 \cos (r)e^{-0.15r}$
    (a): variation of $\mathcal{V}_{\max}(t)$ and
    $\mathcal{V}_{\min}(t)$ with time from $t=0$ to $t=80$. (b): the
    difference $\mathcal{V}_{\max}(t)-\mathcal{V}_{\min}(t)$ tends to
    zero as $t\to\infty$.}
  \label{cos_srt_min_max_v}
\end{figure}

Now we show that the normal velocity $M$ and the gradient
$\mathcal{V}$ of the gas density at the shock tend to become constant
as the computational time increases. As before, let us denote by
$M_{\max}(t)$ and $M_{\min}(t)$, the maximum and minimum of $M$ taken
over $(\xi_1,\xi_2)$ at any time $t$. In an analogous manner we
define $\mathcal{V}_{\max}(t)$ and $\mathcal{V}_{\min}(t)$. In
Figure~\ref{cos_srt_min_max}(a)-(b) we plot the distribution of
$M_{\max}(t),M_{\min}(t)$ and $M_{\max}(t)-M_{\min}(t)$ with respect
to time from $t=0$ to $t=80$. It can be seen that both $M_{\max}(t)$
and $M_{\min}(t)$ decay to one with time and as a result, the
difference $M_{\max}(t)-M_{\min}(t)$ tends to zero asymptotically. We
have given the plot of $\mathcal{V}_{\max}(t),\mathcal{V}_{\min}(t)$
and $\mathcal{V}_{\max}(t)-\mathcal{V}_{\min}(t)$ in
Figure~\ref{cos_srt_min_max_v}(a)-(b). From the figure it is clear
that both $\mathcal{V}_{\max}(t)$ and $\mathcal{V}_{\min}(t)$ approach
zero as time increases.


\subsection{Converging Shock Front Initially in the Shape of a
  Circular Cylinder}
\label{subsec:cylinder}
In this test problem we present the results of simulation of a
cylindrically converging shock front. Even though this is basically a
2-D shock propagation problem, we intent to study it with our 3-D
numerical method.

The initial geometry of the front is a portion of a circular
cylinder of radius two units, i.e.\
\begin{equation}
  \Omega_0\colon x_1^2+x_2^2=4, \ -\frac{\pi}{2}\leq x_3 \leq
  \frac{\pi}{2}.
  \label{cylinder_omega0}
\end{equation}
Initially, the ray coordinates $(\xi_1,\xi_2)$ are chosen as
$\xi_1=x_3$ and $\xi_2=\theta$, where $\theta$ is the azimuthal
angle. Therefore, the initial shock front $\Omega_0$, given in
(\ref{cylinder_omega0}), can be expressed in a parametric form
\begin{equation}
  \Omega_0\colon x_1=2\cos\xi_2, \ x_2=2\sin\xi_2, \ x_3=\xi_1, \
  -\frac{\pi}{2}\leq \xi_1 \leq \frac{\pi}{2}, \ 0 \leq \xi_2 \leq
  2\pi.
  \label{cylinder_ics}
\end{equation}
We have imposed periodic boundary conditions at $\xi_2=0$ and
$\xi_2=2\pi$ and extrapolation boundary conditions at $\xi_1=\pm\pi/2$.

As a result of the particular choice of the ray coordinates
$(\xi_1,\xi_2)$, the unit normal to $\Omega_0$ given by
$\uu{n}_0=(-\cos\xi_2,\sin\xi_2,0)$, points inward and hence the front
converges. If the initial velocity is given a uniform distribution on
$\Omega_0$ as in the previous problems, the front $\Omega_t$ at any
successive time $t$ will remain a circular cylinder with no
interesting geometrical features. Our aim here is to study the
stability of a focusing shock front to perturbations. Hence, the
initial distribution of the normal velocity $M$ is given as a small
perturbation of a constant value, i.e.\
\begin{equation}
  M_0(\xi_1,\xi_2)=1.2+\alpha\cos(\nu\xi_2)
  \label{cylinder_m0}
\end{equation}
with $\alpha=0.05$ and $\nu=8$.

\begin{figure}[htbp]
  \begin{center}
    \includegraphics[height=0.265\textheight]{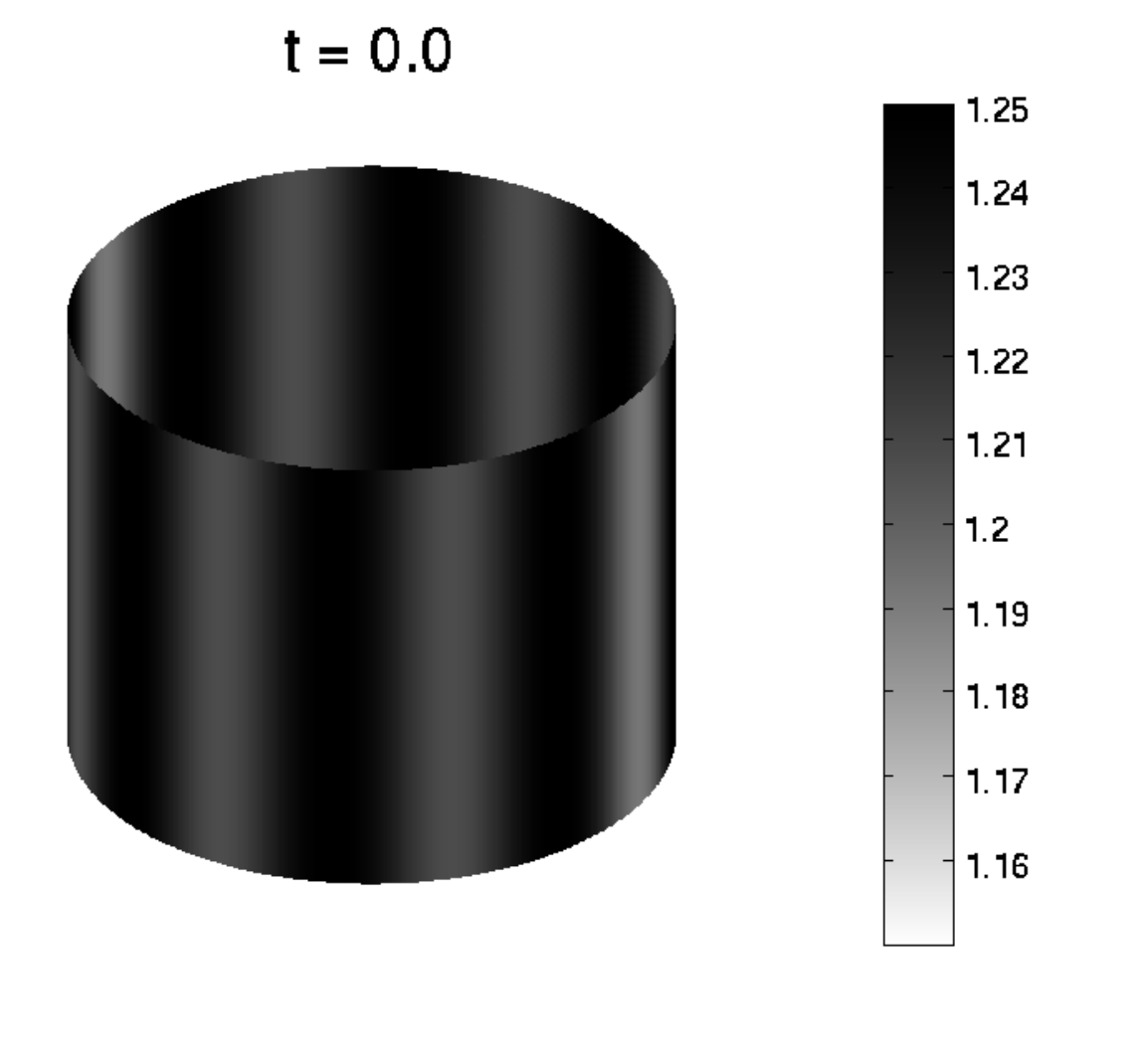}
    \qquad
    \includegraphics[height=0.265\textheight]{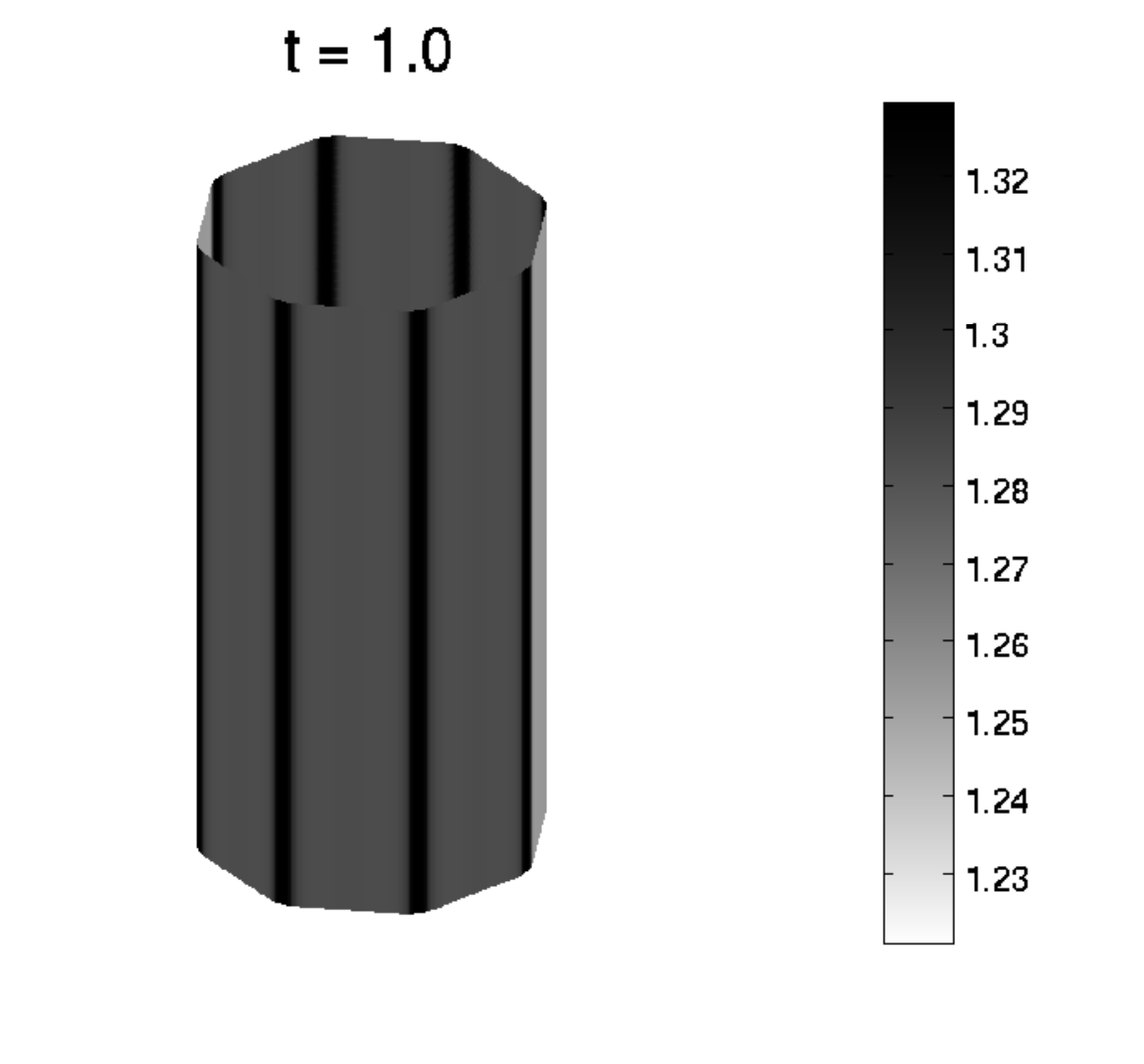}
  \end{center}
  \caption{Cylindrically converging shock front. On the left the
    initial front and on the right the front at time $t=1.0$. The
    grayscale-bar on the right hand side indicates the intensity of
    the normal velocity $M$.}
  \label{cylinder_srt_3d}
\end{figure}

In Figure~\ref{cylinder_srt_3d} we give the plots of the initial shock
front $\Omega_0$ and the shock front $\Omega_t$ at time $t=1.0$. The
shock fronts are coloured using the variation of the normal velocity
$M$, with the grayscale-bar on the right indicating the values of
$M$. Note that the initial normal velocity $M_0$ has a periodic
variation with a maximum value $1.25$ and a minimum value
$1.05$, cf.\ \eqref{cylinder_m0}. Those portions of the front where
$M_0$ has maximum value moves inwards faster and it results in a
distortion of the circular shape of $\Omega_0$. From the shock front
at $t=1$ in Figure~\ref{cylinder_srt_3d}, it can be observed that
sixteen fully developed vertical kink lines are formed on the shock
front and kink lines have not yet interacted. Clearly, the number of
plane sides and kink lines on $\Omega_t$ at various times will depend
on the value of parameter $\nu$ in \eqref{cylinder_m0}. The shock
front assumes the shape of a polygonal cylinder and as time
progresses the Mach stem-like surfaces on the shock front interacts,
with the formation of new kink lines and this process repeats. It
has to be remarked that our numerical results are well in accordance
with the experimental results of \cite{takayama-etal}, where the
authors have reported that a converging shock front with a small
perturbation assumes the shape of a polygon. This test problem, hence,
gives another yet evidence for the efficacy of SRT to produce
physically realistic geometrical features.

\begin{figure}
  \begin{center}
    \includegraphics[height=0.4\textheight]{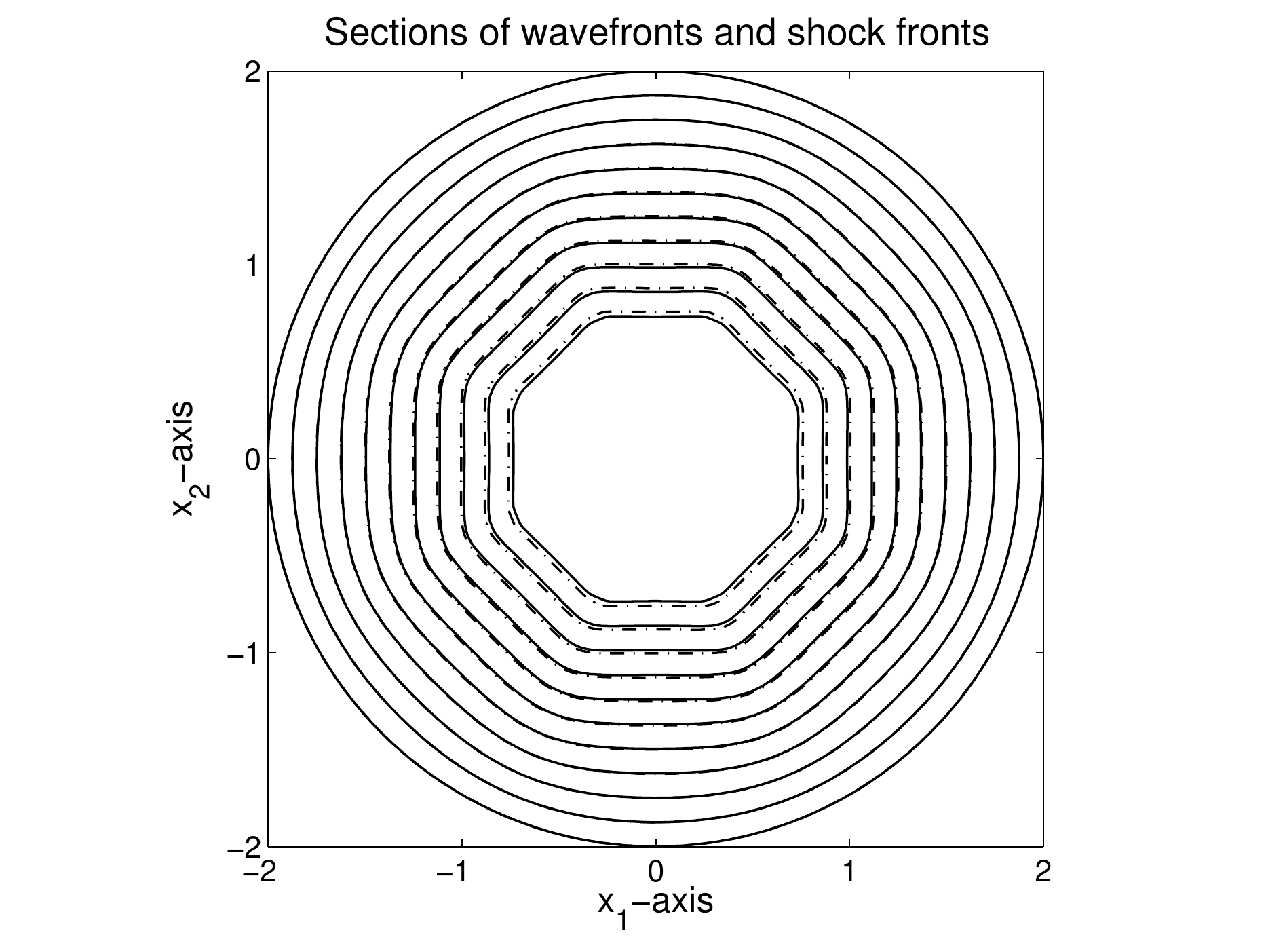}
  \end{center}
  \caption{Cross-sections of the shock fronts (broken lines) and the
    nonlinear wavefront (continuous lines) at $t=0,0.1,0.2,\ldots,1.0$
    by $x_3=0$ plane}
  \label{cylinder_srt_wnlrt_sect}
\end{figure}

An examination of the shock front at $t=1.0$ in
Figure~\ref{cylinder_srt_3d} and the corresponding nonlinear
wavefront in \cite{arun-ct} shows that both the fronts have
analogous geometry. We also refer the reader to \cite{arun-ct} for
more details on the transient geometries between the initial and
final configurations and a quantification of the focusing process.
Nevertheless, in order to point out the important difference in
the case of a shock front, in Figure~\ref{cylinder_srt_wnlrt_sect}
we give the successive cross-sections (say by $x_3=0$ plane) of
the shock fronts and the corresponding nonlinear wavefronts at
times $t=0,0.1,0.2,\ldots,1.0$. Note that initially the normal
speeds $m$ and $M$ of both the fronts have the same initial
values. These speeds increase due to convergence of the fronts,
however, in the case of a shock front there is also a decay in the
shock speed due to interaction of the shock with the nonlinear
wavefronts, cf.\ ({\it also}) the second term in
\eqref{2d_srt_cons3}. Therefore, as time increases, the shock
front lags behind the nonlinear wavefront as seen in
Figure~\ref{cylinder_srt_wnlrt_sect}.

\subsection{Spherically Converging Shock Front}
\label{subsec:sphere}

In this problem we consider the propagation of a spherically
converging shock front. The initial geometry of the shock front is a
sphere of radius two units. The ray coordinates are chosen to be
$\xi_1=\pi-\phi$ and $\xi_2=\theta$, where $\theta$ is the azimuthal
angle and $\phi$ is the polar angle. Therefore, the parametric
representation of the initial shock front $\Omega_0$ is
\begin{equation}
  x_1=2\sin\xi_1\cos\xi_2,\ x_2=2\sin\xi_1\sin\xi_2, \ x_3=-2\cos\xi_1.
\end{equation}
In order to avoid the singularities at $\phi=0$ and $\phi=\pi$,  we
remove these points. Therefore, our computational domain is
$[\pi/15,14\pi/15]\times[0,2\pi]$. As in the previous problem, we
choose the initial velocity distribution as a small perturbation of a
constant value, i.e.\
\begin{equation}
  M_0(\xi_1,\xi_2)=1.2+\alpha\cos(\nu_1\xi_1)\cos(\nu_2\xi_2)
  \label{sphere_m0}
\end{equation}
with $\alpha=0.05,\nu_1=4,\nu_2=8$.

\begin{figure}[htbp]
  \begin{center}
    \begin{tabular}{cc}
      \includegraphics[width=0.45\textwidth]{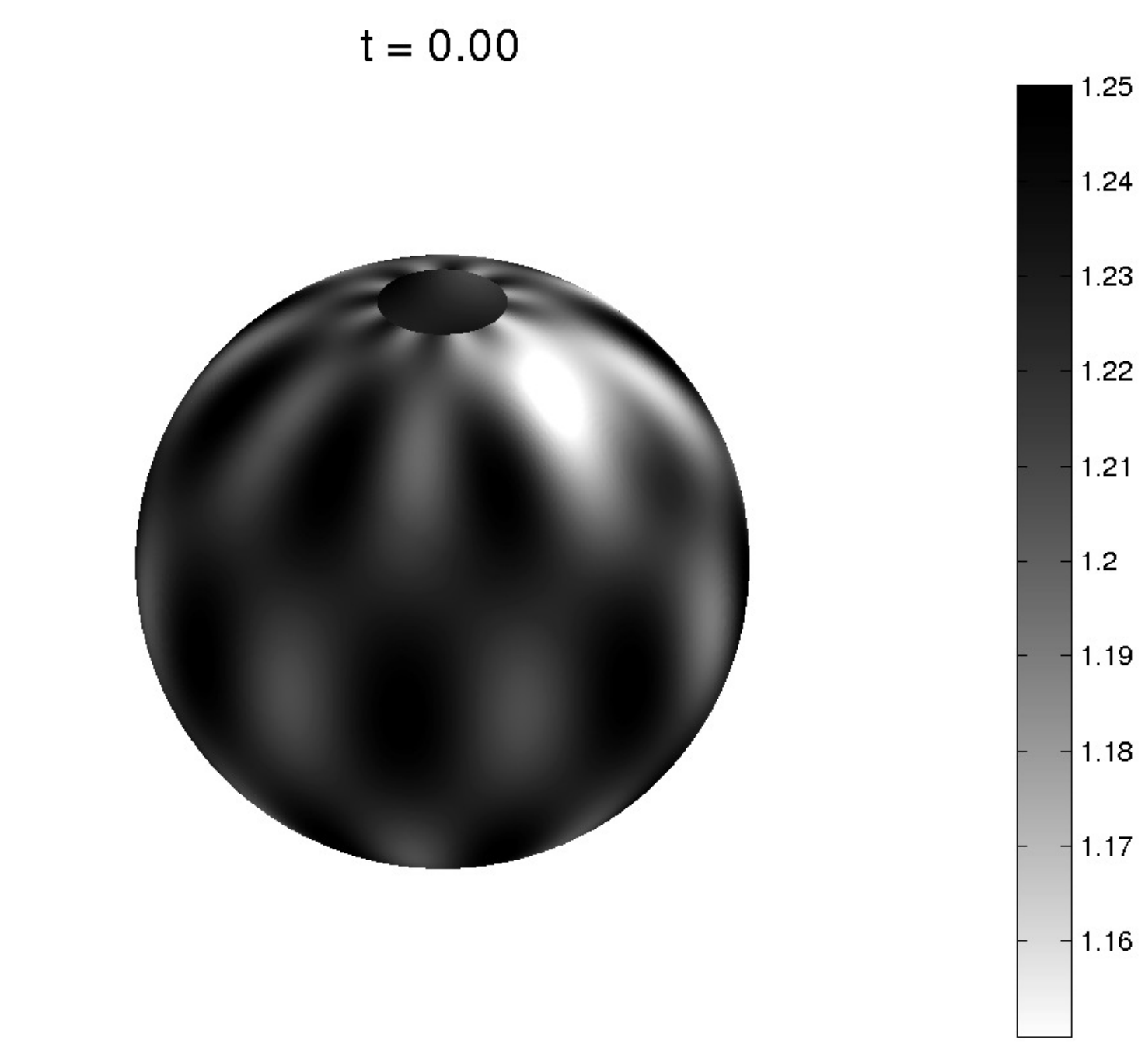}&
      \includegraphics[width=0.45\textwidth]{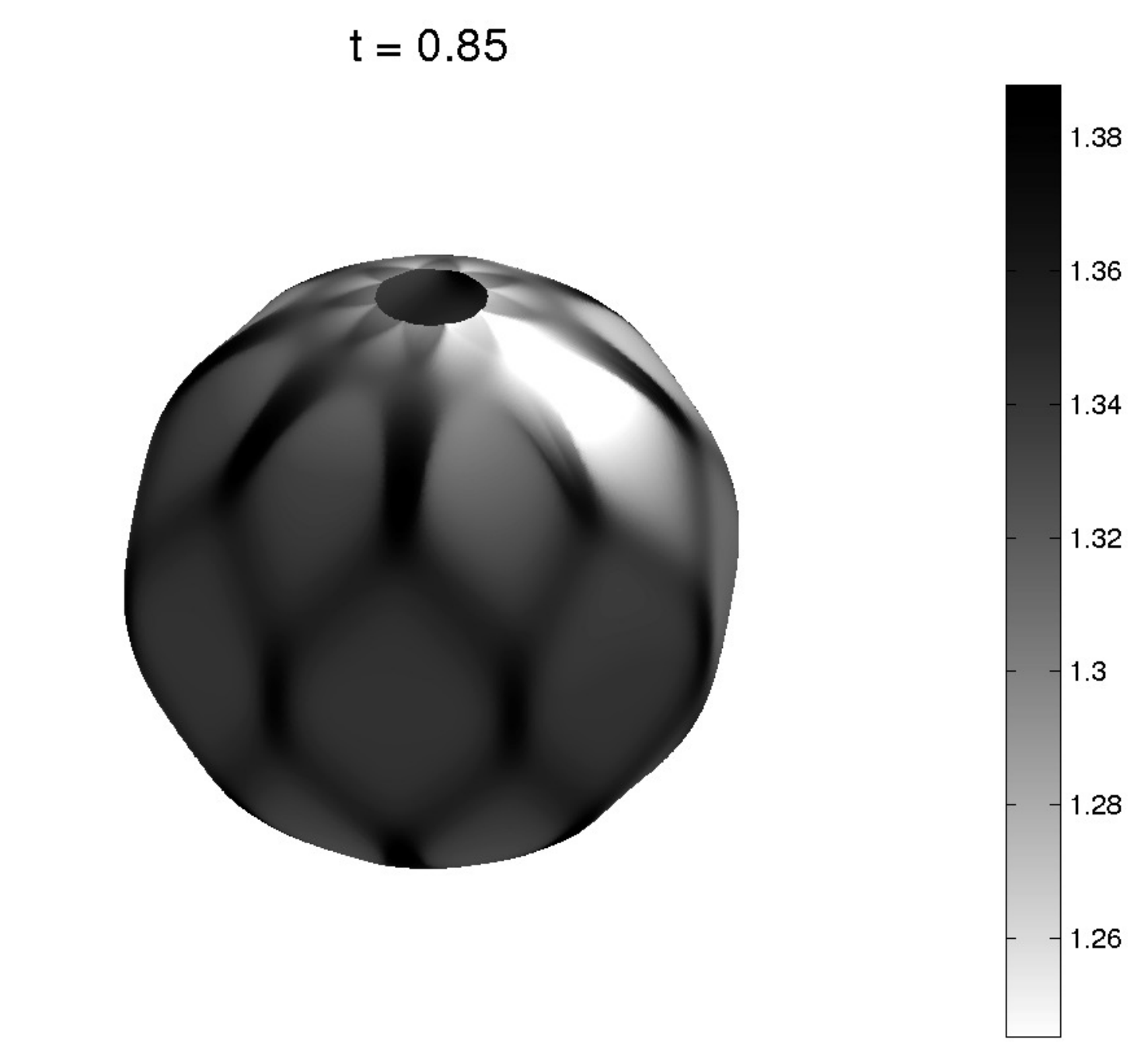}
    \end{tabular}
  \end{center}
  \caption{Spherically converging shock front at $t=0.85$. The
    grayscale bar on the right represents the distribution of $M$.}
  \label{sphere_srt_3d}
\end{figure}

The 3-D plots of the initial shock front and the one at time $t=0.85$
are given in Figure~\ref{sphere_srt_3d}. The shock fronts are coloured
using the variation of the normal velocity $M$, with the grayscale-bar on
the right indicating the values of $M$. It can be observed from the
figure that as the front starts focusing, it develops several kink
curves and its spherical shape gets distorted, with the formation of
facets. The final shape of the shock front is almost a polyhedron.

In this test problem also we have a observed qualitatively similar
behaviour of a shock front and nonlinear wavefront. The major
important difference is that the shock front slightly lags behind
the nonlinear wave due to the effect on the shock of the flow behind
the shock, as explained in the case of a cylindrically converging
shock.
\begin{remark}
  It has to be emphasised that we can continue the computations
  further, however, the results would not be physically realistic as
  the weak shock ray theory breaks down for larger values of $M-1$.
\end{remark}
\begin{remark}
  As in the case nonlinear wavefront in \cite{arun-ct}, we notice that
  cylindrically and spherically converging shocks take the polygonal and
  polyhedral shapes, respectively. These two shapes are stable
  configurations for the two cases. The stability here does not mean
  that once one of these configuration is formed, the shock front
  remains in this configuration for all time. It may change again into
  another such configuration.
\end{remark}

\subsection{Some More Comparison Between 3-D WNLRT and 3-D SRT}
\label{subsec:comp}

By comparing the results of numerical experiments presented above with
those of 3-D WNLRT reported in (\cite{arun-ct}), we infer that the
geometrical shapes of nonlinear wavefronts and shock fronts are more
or less qualitatively similar. However, looking at these shapes alone,
it is not possible say much about the ultimate results as
$t\to\infty$.

When the initial shape is periodic in $x_1$ and $x_2$ given by
\ref{sine_omega0}, we have examined the behaviour of a nonlinear
wavefront in \cite{arun-ct} and that of a shock front in
subsection~\ref{subsec:corrug} for $t\to\infty$ and found that both
are corrugationally stable. This means that both these fronts tend to
become planar after a long time, which in turn shows that the mean
curvature ${\it \Omega}$ approaches zero. Writing the differential
form of \eqref{3d_wnlrt_cons3} and letting ${\it \Omega}\to 0$, i.e.\
the mean curvature approaching zero yields
\begin{equation}
  m(\xi_1,\xi_2,t)\to \const, \ \mbox{as} \ t\to\infty.
  \label{m_limit}
\end{equation}
In order that the nonlinear wavefront be corrugation stable, this
constant must be the same along all rays. It has been observed in
\cite{arun-ct} that even though the maximum and minimum of the
front velocity $m$ of a nonlinear wavefront heavily oscillate about
their initial values, they finally approach their initial mean values
as $t\to\infty$. It is, therefore, very interesting to note that the
numerical computation suggests the constant in (\ref{m_limit}) is
$m_0$ for the case when $m_0$ is constant everywhere on the initial
wavefront $\Omega_0$.

We now proceed to investigate the long term of behaviour of
the perturbations on a 3-D corrugationally stable plane shock
front. First, we note that the numerical results show that both
$M_{\max}$ and $M_{\min}$ decrease to one with increasing
time. Therefore, $M\to 1$ as $t\to \infty$, cf.\
subsections~\ref{subsec:corrug} and \ref{subsec:axisym}. The results
also show that the gradient $\mathcal{V}$ of the pressure decays to
zero as $t\to\infty$. Analogous results were reported also in the 2-D
case in \cite{monica}. Thus, we can see a major difference in the
decay of the shock amplitude from the corresponding results of 3-D
WNLRT.

\begin{figure}
  \begin{center}
    \includegraphics[width=0.6\textwidth]{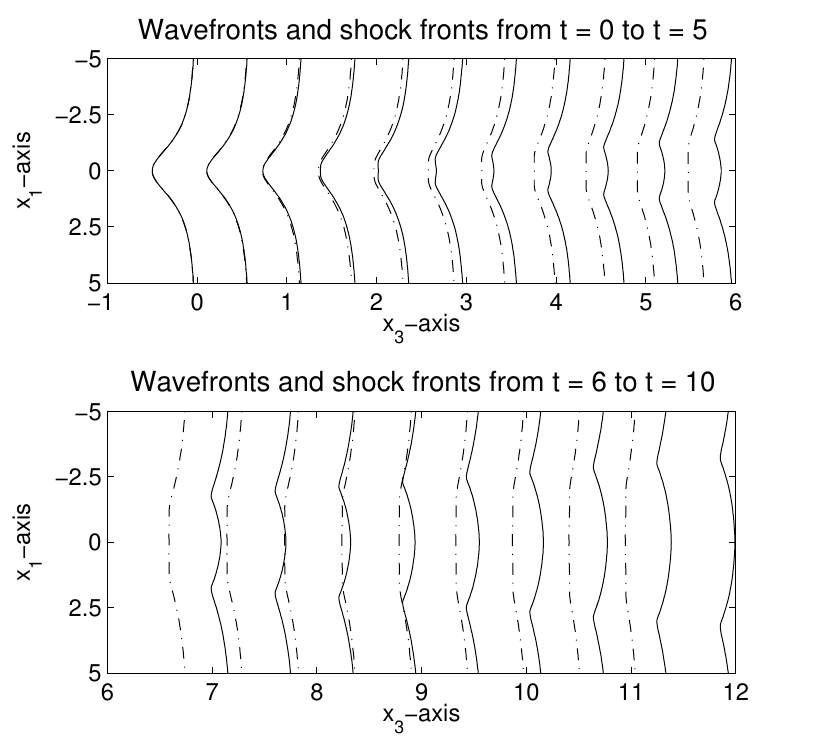}
  \end{center}
  \caption{Comparison of the results by 3-D WNLRT and 3-D SRT. Figure
    on the top: from $t=0.0$ to $t=5.0$ and bottom: from $t=6.0$ to
    $t=10.0$. The solid lines represent successive positions of the
    nonlinear wavefront obtained by 3-D WNLRT and dotted lines are
    those of the shock front from 3-D SRT. The kinks can be noticed on
    both the fronts. The nonlinear wavefront overtakes the shock and
    its central portion bulges out so that the kinks on it are
    very prominent.}
  \label{wnlrt_srt_comp}
\end{figure}

In order to do some more comparison of results for  a nonlinear
wavefront and shock front, we choose the initial geometry of the
nonlinear wavefront and shock front to be an axisymmetric dip given
in (\ref{dip_omega0}) with $\kappa=0.5,\alpha=\beta=1.5$. We take the
same initial values of $M$ and $m$, both equal to 1.2. This choice is
different from that in \cite{monica}, where $\mu$ defined by
\eqref{mu_defn} was chosen to be the same. Note that $M$ is given by
\eqref{M_V}, whereas for a nonlinear wavefront
\begin{equation}
  \label{eq:m_mu}
  m:=1+\varepsilon\frac{\gamma+1}{2}\mu.
\end{equation}
The computations are done with 3-D WNLRT and 3-D SRT up to a time
$t=10$. A comparison of the results obtained is presented in
Figure~\ref{wnlrt_srt_comp}, where we have plotted the successive
wavefronts and shock fronts in the section $x_2=0$ from $t=0$ to
$t=10$ in a time step of $0.5$. In the figure, the solid lines
represent the successive nonlinear wavefronts and dotted lines are the
shock fronts. The figure clearly shows that from time $t=2.0$ onwards,
the nonlinear wavefront overtakes the corresponding shock front. We
also notice that the central portion of the nonlinear wavefront bulges
out and the two kinks move apart faster than those on the shock
front.

\begin{figure}
  \begin{center}
    \begin{tabular}{cc}
      \includegraphics[width=0.45\textwidth]{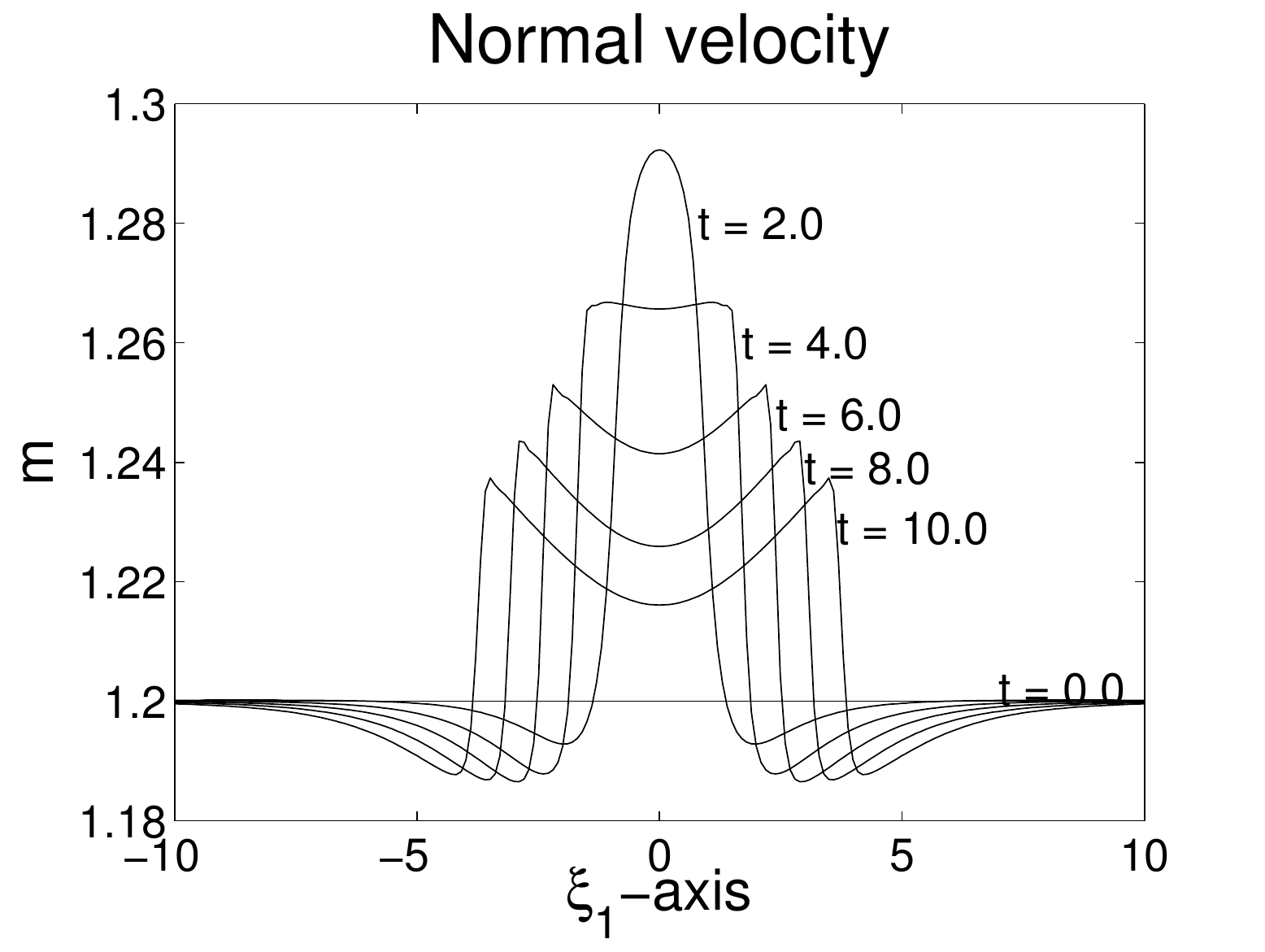}&
      \includegraphics[width=0.45\textwidth]{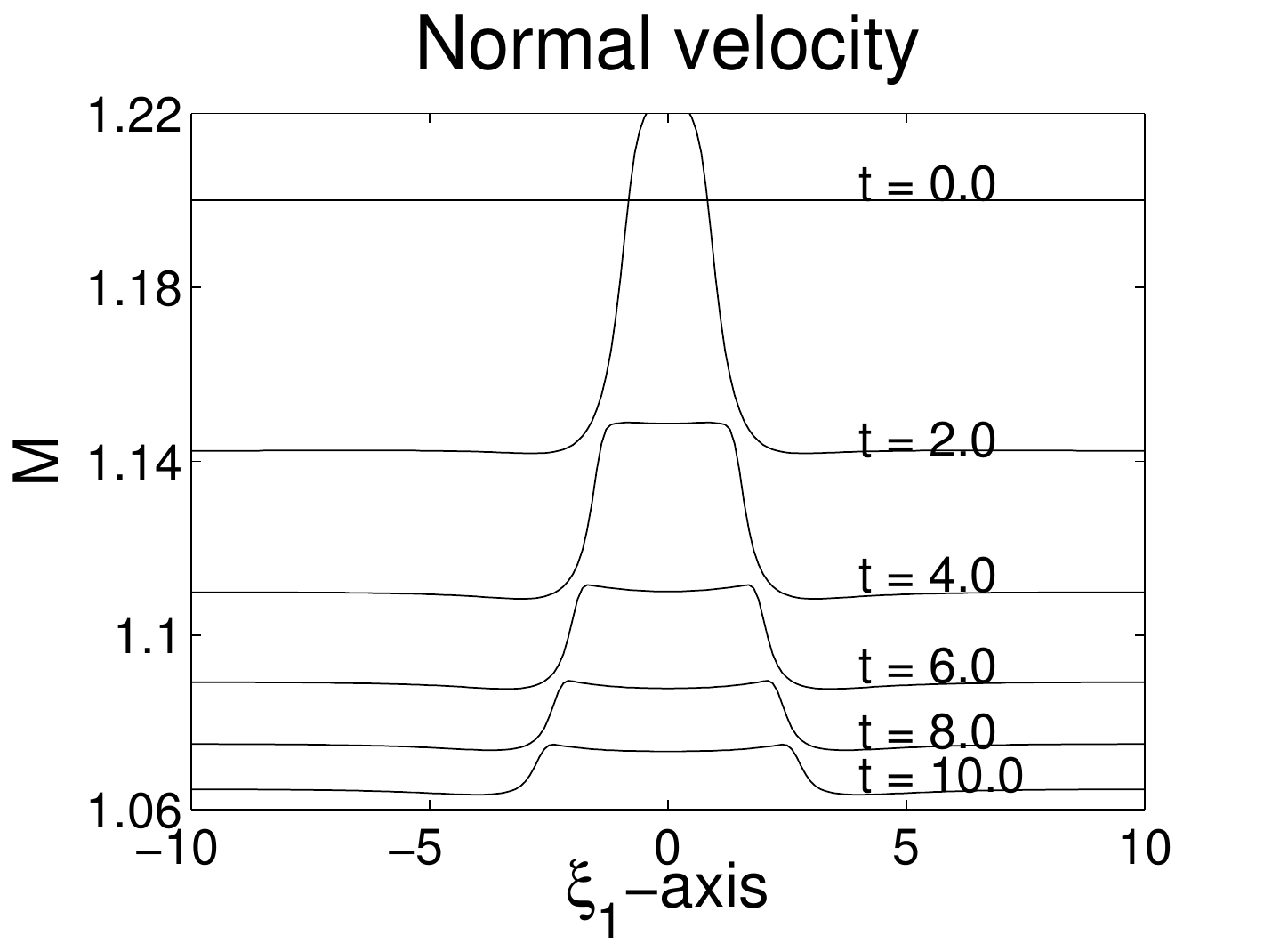}\\
      (a) & (b)
    \end{tabular}
  \end{center}
  \caption{Comparison of the Mach number distribution on the fronts at
    times $0,2,4,6,8,10$ with the same amplitude distribution on the
    initial fronts. (a): results obtained by 3-D WNLRT. (b): those of
    3-D SRT.}
  \label{wnlrt_srt_m_comp}
\end{figure}

We also present in Figure~\ref{wnlrt_srt_m_comp} the graphs of the
normal velocity $m$ of the nonlinear wavefront and $M$ of that of the
shock front. In the figure we have plotted both $m$ and $M$ at times
$0,2,4,6,8,10$ with the same constant initial values on the respective
fronts. The Mach number at the centre of the fronts initially rises
considerably in both the cases but it becomes constant on the central
disc at time $t=4$ for the nonlinear wavefront. As the rays starts
diverging from the bulged central portion,
cf. Figure~\ref{wnlrt_srt_m_comp}(a), we can see that $m$ reduces at
the central portion from $t=6$ onwards. However, as seen in all
previous cases with $\mathcal{V}_0>0$, the shock Mach number (and
hence the shock amplitude) decreases with time on all parts of the
shock front governed by 3-D SRT.

\section{Concluding Remarks}
\label{sect:remarks}

Using the 3-D KCL based SRT we have successfully calculated the
evolution of several shock surfaces, starting from a wide range of
very interesting initial shapes. The proposed theory and numerical
solution procedure correctly takes into account of the effect of the
flow behind the shock, which is very important for weak shocks. It has
been observed that the geometry, position and amplitude of a shock
depends on its initial position, amplitude distribution and gradient of
flow variables behind it. A comparison of the results with those of a
weakly nonlinear wavefront shows that a weak shock front and a weakly
nonlinear wavefront are topologically same.

With the aid of extensive numerical simulations, we have been able to
verify the corrugational stability of a planar shock and have quantified
the decay of small perturbations in the shape of a planar front. The
simulation of cylindrically and spherically converging shock fronts
shows that they assume polygonal and polyhedral shapes,
receptively. These configurations appear to be stable in the sense
that every such periodic distribution of $M$ on an initially
cylindrical and spherical shock is likely to lead to polygonal and
polyhedral shapes, which later on may get transformed to other similar
shapes due to interaction of kink lines. However, our computations was
not be continued further as the small amplitude assumption is violated
because the value of $M$ increased too much due to radial
convergence. Further calculation requires a new formulation of both
WNLRT and SRT for the arbitrary amplitude case, along the lines of
\cite{prasad-parker}, and is a subject matter for future research.

In all the numerical simulations, particularly in long time
computations, we have observed an interesting phenomenon of the
persistence of kinks curves. As we have seen in all the cases
considered by us, a kink curve may appear on an initially smooth shock
front, but once it is formed it persists till it meets another
kink curve. The persistence of a kink follows from the similar
property of a shock in a genuinely nonlinear characteristic field; a
shock, once formed cannot terminate at a finite distance in the
$(\xi_1,\xi_2,t)$-space; see \cite{prasad-book}. We further notice
that interaction of a pair of kink curves (whether they correspond to
the same characteristic field or different characteristic fields)
always produces another pair of kink curves; see
\cite{baskar-riemann}.



\appendix
\section{Jacobian Matrices}
\label{sect:append_a}
A quasilinear form of the balance equations \eqref{3d_srt_cons} of 3-D
SRT can be written in a matrix form
\begin{equation}
  \tilde{A}V_t+\tilde{B}^{(1)}V_{\xi_1}+\tilde{B}^{(2)}V_{\xi_2}=\tilde{C},
\end{equation}
where $V=(U_1,U_2,V_1,V_2,M,G_1,G_2,\mathcal{V})^T$. The Jacobian
matrices
$\tilde{A}=(\tilde{a}_{ij}),\tilde{B}^{(1)}=({\tilde{b}}^{(1)}_{ij})$
and $\tilde{B}^{(2)}=({\tilde{b}}^{(2)}_{ij})$ are of size $8\times 8$
and the vector $\tilde{C}$ belongs to $\mathbb{R}^8$. The nonzero
elements of $\tilde{A},\tilde{B}^{(1)},\tilde{B}^{(2)}$ and
$\tilde{C}$ are given below.
\begin{gather*}
  {\tilde{a}}_{11}={\tilde{a}}_{22}=G_1, \
  {\tilde{a}}_{33}=G_2, \ {\tilde{a}}_{44}=G_2, \
  {\tilde{a}}_{51}=-\frac{1}{U_3}G_1G_2N_2\cot\chi, \
  {\tilde{a}}_{52}=\frac{1}{U_3}G_1G_2N_1\cot\chi, \\
  {\tilde{a}}_{53}=\frac{1}{V_3}G_1 G_2N_2\cot\chi, \
  {\tilde{a}}_{54}=-\frac{1}{V_3}G_1G_2N_1\cot\chi, \
  {\tilde{a}}_{55}=\frac{2M}{M-1}G_1G_2, \ {\tilde{a}}_{56}=G_2, \
  {\tilde{a}}_{57}=G_1, \\
  {\tilde{a}}_{66}={\tilde{a}}_{77}=1,\
  {\tilde{a}}_{81}=-\frac{1}{U_3}G_1G_2N_2\cot\chi, \
  {\tilde{a}}_{82}=\frac{1}{U_3}G_1G_2N_1\cot\chi, \
  {\tilde{a}}_{83}=\frac{1}{V_3}G_1 G_2N_2\cot\chi, \\
  {\tilde{a}}_{84}=-\frac{1}{V_3}G_1G_2N_1\cot\chi, \
  {\tilde{a}}_{85}=2G_1G_2,\ {\tilde{a}}_{86}=G_2,\ {\tilde{a}}_{87}=G_1,\
  {\tilde{a}}_{88}=\frac{2G_1G_2}{\mathcal{V}}
\end{gather*}
\begin{gather*}
  {\tilde{b}}^{(1)}_{11}=-\frac{M}{U_3}(U_1U_2+N_1N_2)\cot\chi, \
  {\tilde{b}}^{(1)}_{12}=\frac{M}{U_3}(U^2_1+N^2_1-1)\cot\chi, \\
  {\tilde{b}}^{(1)}_{13}=\frac{M}{V_3\sin\chi}(U_2 V_1+N_1N_2\cos\chi), \
  {\tilde{b}}^{(1)}_{14}=-\frac{M}{V_3\sin\chi}(U_1V_1+(N^2_1-1)\cos\chi),\
  {\tilde{b}}^{(1)}_{15}=-N_1.
\end{gather*}
\begin{gather*}
  {\tilde{b}}^{(1)}_{21}=-\frac{M}{U_3}(U^2_2+N^2_2-1)\cot\chi, \
  {\tilde{b}}^{(1)}_{22}=\frac{M}{U_3}(U_1U_2+N_1N_2)\cot\chi, \\
  {\tilde{b}}^{(1)}_{23}=\frac{M}{V_3\sin\chi}(U_2V_2+(N_2^2-1)\cos\chi), \
  {\tilde{b}}^{(1)}_{24}=-\frac{M}{V_3\sin\chi}(U_1V_2+N_1N_2\cos\chi), \
  {\tilde{b}}^{(1)}_{25}=-N_2.
\end{gather*}
\begin{gather*}
  {\tilde{b}}^{(1)}_{61}=-\frac{M}{U_3\sin\chi}(V_2-U_2\cos\chi), \
  {\tilde{b}}^{(1)}_{62}=\frac{M}{U_3\sin\chi}(V_1-U_1\cos\chi).
\end{gather*}
\begin{gather*}
  {\tilde{b}}^{(2)}_{31}=-\frac{M}{U_3\sin\chi}(U_1V_2+N_1N_2\cos\chi), \
  {\tilde{b}}^{(2)}_{32}=\frac{M}{U_3\sin\chi}(U_1V_1+(N^2_1-1)\cos\chi),\\
  {\tilde{b}}^{(2)}_{33}=\frac{M}{V_3}(V_1V_2+N_1N_2)\cot\chi, \
  {\tilde{b}}^{(2)}_{34}=-\frac{M}{V_3}(V^2_1+N^2_1-1)\cot\chi, \
  {\tilde{b}}^{(2)}_{35}=-N_1.
\end{gather*}
\begin{gather*}
  {\tilde{b}}^{(2)}_{41}=-\frac{M}{U_3\sin\chi}(U_2V_2+(N^2_2 -1)\cos\chi), \
  {\tilde{b}}^{(2)}_{42}=\frac{M}{U_3\sin\chi}(U_2V_1+N_1N_2\cos\chi),\\
  {\tilde{b}}^{(2)}_{43}=\frac{M}{V_3}(V^2_2+N^2_2-1)\cot\chi, \
  {\tilde{b}}^{(2)}_{44}=-\frac{M}{V_3}(V_1V_2+N_1N_2)\cot\chi, \
  {\tilde{b}}^2_{45}=-N_2.
\end{gather*}
\begin{gather*}
  {\tilde{b}}^{(2)}_{73}=\frac{M}{V_3\sin\chi}(U_2-V_2\cos\chi), \
  {\tilde{b}}^{(2)}_{74}=-\frac{M}{V_3\sin\chi}(U_1-V_1\cos\chi).
\end{gather*}
\begin{gather*}
  \tilde{c}_7=-\frac{2M}{M-1}G_1G_2\mathcal{V}, \
  \tilde{c}_8=-2MG_1G_2\mathcal{V}.
\end{gather*}

\section{Some Explanations of the KCL Based SRT}
\label{sect:letter b}
1. There are two parts in our method. The first part (purely
geometrical) consists of KCL, in which the six jump relations imply
conservation of distances in $ x_1, x_2$ and $x_3$ directions (and
hence in any arbitrary direction in {\bf x}-space) across a shock, see
theorem 3.1 on page 297 of \cite{arun-prasad-wave}, also the section
3.3.3 of \cite{prasad-book} for a detailed discussion. 

The second part consists of the two closure relations
\eqref{3d_srt_mu1}-\eqref{3d_srt_mu2} or \eqref{3d_srt_cons3} and
\eqref{3d_srt_cons4}, which are derived from the Euler equations of
gas-dynamics. Out of these two closure relations, the first one
explicitly represents conservation of energy along a shock ray. We
note that the 4 jump relations of a 2-D shock (representing
conservation of mass, momentum and energy) express all quantities
behind the shock in terms of the unit normal $\bf N$ to the shock and
a single unknown, say the Mach number $M$ of the shock. Time rate of
change along a ray of $\bf N$ is a geometrical relation and is implied
by KCL (see section 5, \cite{arun-prasad-wave}) . Therefore, it is
sufficient know the rate of change of $M$ along a shock ray. This is
the basis of the derivation of the infinite system of equations by
Grinfel'd \cite{grinfeld} and Maslov \cite{maslov}. Thus the two
equations for $M$ and $\mathcal{V}$ along a shock ray in NTSD,
implicitly take into account of conservation of mass and momentum
also. 

The geometrical conservation laws give the jump relations across a
kink. The Euler equations, which are used to deduce the two
compatibility conditions in conservation form, cannot account for the
kink phenomenon. Initially, we had only the transport equation for $M$
along a shock ray in differential form [28]. Thus, we started with
Euler equations. It was the need for capturing the kinks (geometrical
singularities) on a shock that led to the formulation (or discovery)
of 2-D KCL \cite{morton} in 1992 and 3-D KCL \cite{giles} in 1995. 

2. At the end of the section 2 on page 5 we have stated "For a
discussion on the validity of NTSD and its application to 2-D problems
we refer the reader to
\cite{baskar-jfm,kevlahan,monica,prasad-book}. It is interesting to
note that the NTSD gives quite good results even for a shock of
arbitrary strength, which has been verified for a 1-D piston problem
in \cite{lazarev-prasad-sing}. In addition, it has also been observed
in \cite{lazarev-prasad-sing} that the NTSD takes less than $0.5\%$ of
the computational time needed by a typical finite difference method
applied to the Euler equations." Let me elaborate some of these below.  

Before we found KCL, we had the NTSD in differential form and that is
equivalent to KCL based shock ray theory for smooth shocks. While
formulating the NTSD with a model equation, we showed a good agreement
with exact solution \cite{ravindran-prasad} and a companion paper in
the same journal. Ravindran, Sunder and Prasad also found good
agreement in 1994 (see reference to this work in \cite{kevlahan}). The
first detailed attempt to verify and compare our method with other
methods was done by Kevlahan \cite{kevlahan}. Though he used the NTSD
in differential form, he had a method to locate the kinks (see page
180 in \cite{kevlahan}) and this method indeed gives the kinks
location (see \cite{prasad-book}, page 126). Kevlahan cocludes "{\it
  The theory is tested against known analytical solutions for
  cylindrical and plane shocks, and against a full direct numerical
  simulation (DNS) of a shock propagating into a sinusoidal shear
  flow. The test against DNS shows that the present theory accurately
  predicts the evolution of a moderately weak shock front, including
  the formation of shock-shocks due to shock focusing. The theory is
  then applied to the focusing of an initially parabolic shock and he
  finds that the shock shapes (by this theory) agreed well with the
  experimental results}". 

Since our theory is valid for a weak shock, we cannot compare with
similarity solution of Guderley but, Kevlahan even compares a higher
order theory of ours and finds remarkably good agreement.  

A later reference \cite{baskar-jfm} using 2-D KCL contains another
extensive comparison of results with numerical solution of Euler
equations. 

3. The choice that the characteristic length $L$ is of the order of
the distance over which 3-D SRT is valid, though it may look a bit vague
but is an appropriate choice. A clear choice could have been the
inverse of the mean curvature of the initial geometry of the
shock. But due to appearance of the kinks, this information is lost in
examples in subsection~\ref{subsec:corrug} and
subsection~\ref{subsec:axisym} and solution is valid for almost
infinite time.  In the example in subsection~\ref{subsec:non_axisym}
our result is calculated up to time $t=10$ and shock has moved
approximately by a distance 11.08. However, we stop the computation at
this stage as the value of $M$ increases beyond the validity of weak
shock assumption. In the examples in subsection~\ref{subsec:cylinder}
and subsection~\ref{subsec:sphere}, we could have taken $L$ to be the
radius of the cylinder and that of the sphere respectively - these are
of the order of the distance traveled by the shock. Computation in
these cases were stopped because the value of $M-1$ became too large
for the weak shock assumption to be valid.

\bibliography{references}

\begin{thebibliography}{10}

\bibitem{anile-russo-1986}
A.~M. Anile and G.~Russo.
\newblock Corrugation stability for plane relativistic shock waves.
\newblock {\em Phys. Fluids}, 29:2847--2852, 1986.

\bibitem{arun-ct}
K.~R. Arun.
\newblock A numerical scheme for three-dimensional front propagation and
  control of {Jordan} mode.
\newblock {\em SIAM J. Sci. Comput.}, 34:B148--B178, 2012.

\bibitem{arun-etal-siam}
K.~R. Arun, M.~{Luk\'a\v{c}ov\'a}-{Medvi\softd ov\'a}, P.~Prasad, and S.~V.
  {Raghuram}a Rao.
\newblock An application of {3-D} kinematical conservation laws: propagation of
  a three dimensional wavefront.
\newblock {\em SIAM J. Appl. Math.}, 70:2604--2626, 2010.

\bibitem{arun-prasad-wave}
K.~R. Arun and P.~Prasad.
\newblock {3-D} kinematical conservation laws ({KCL}): evolution of a surface
  in $\mathbb{R}^3$-in particular propagation of a nonlinear wavefront.
\newblock {\em Wave Motion}, 46:293--311, 2009.

\bibitem{arun-prasad-eig}
K.~R. Arun and P.~Prasad.
\newblock Eigenvalues of kinematical conservation laws ({KCL}) based 3-{D}
  weakly nonlinear ray theory ({WNLRT}).
\newblock {\em Appl. Math. comput.}, 217:2285--2288, 2010.

\bibitem{baskar-riemann}
S.~Baskar and P.~Prasad.
\newblock Riemann problem for kinematical conservation laws and geometrical
  features of nonlinear wavefronts.
\newblock {\em IMA J. Appl. Math.}, 69:391--420, 2004.

\bibitem{baskar-jfm}
S.~Baskar and P.~Prasad.
\newblock Propagation of curved shock fronts using shock ray theory and
  comparison with other theories.
\newblock {\em J. Fluid Mech.}, 523:171--198, 2005.

\bibitem{baskar-sonic}
S.~Baskar and P.~Prasad.
\newblock Formulation of the sonic boom problem by a maneuvering aerofoil as a
  one parameter family of cauchy problems.
\newblock {\em Proc. Indian Acad. Sci. Math. Sci.}, 116:97--119, 2006.

\bibitem{choquet-bruhat}
Y.~Choquet-Bruhat.
\newblock Ondes asymptotiques et approch\'ees pour des syst\`emes d'\'equations
  aux d\'eriv\'es partielles non lin\'eaires.
\newblock {\em J. Math. Pures. Appl.}, 48:117--158, 1969.

\bibitem{evans-hawley}
C.~R. Evans and J.~F. Hawley.
\newblock Simulation of magnetohydrodynamic flows-{A} constrained transport
  method.
\newblock {\em Astrophys. J.}, 332:659--677, 1988.

\bibitem{giles}
M.~B. Giles, P.~Prasad, and R.~Ravindran.
\newblock Conservation forms of equations of three dimensional front
  propagation.
\newblock Technical report, Department of Mathematics, Indian Institute of
  Science, Bangalore, India, 1995.

\bibitem{grinfeld}
M.~A. Grinfel'd.
\newblock Ray method for calculating the wavefront intensity in nonlinear
  elastic material.
\newblock {\em PMM J. Appl. Math. Mech.}, 42:958--977, 1978.

\bibitem{jiang-shu}
G.-S. Jiang and C.-W. Shu.
\newblock Efficient implementation of weighted {ENO} schemes.
\newblock {\em J. Comput. Phys.}, 126:202--228, 1996.

\bibitem{kevlahan}
N.~K.-R. Kevlahan.
\newblock The propagation of weak shocks in non-uniform flows.
\newblock {\em J. Fluid Mech.}, 327:161--197, 1996.

\bibitem{kurganov-tadmor}
A.~Kurganov and E.~Tadmor.
\newblock New high-resolution central schemes for nonlinear conservation laws
  and convection-diffusion equations.
\newblock {\em J. Comput. Phys.}, 160:241--282, 2000.

\bibitem{lazarev-prasad-sing}
M.~P. Lazarev, P.~Prasad, and S.~K. Sing.
\newblock An approximate solution of one-dimensional piston problem.
\newblock {\em Z. Angew. Math. Phys.}, 46:752--771, 1995.

\bibitem{maslov}
V.~P. Maslov.
\newblock Propagation of shock waves in an isentropic non-viscous gas.
\newblock {\em J. Sov. Math.}, 13:119--163, 1980.

\bibitem{monica}
A.~Monica and P.~Prasad.
\newblock Propagation of a curved weak shock.
\newblock {\em J. Fluid Mech.}, 434:119--151, 2001.

\bibitem{morton}
K.~W. Morton, P.~Prasad, and R.~Ravindran.
\newblock Conservation forms of nonlinear ray equations.
\newblock Technical report, Department of Mathematics, Indian Institute of
  Science, Bangalore, India, 1992.

\bibitem{prasad-1975}
P.~Prasad.
\newblock Approximation of perturbation equations of a quasilinear hyperbolic
  system in the neighbourhood of a bicharacteristic.
\newblock {\em J. Math. Anal. Appl.}, 50:470--482, 1975.

\bibitem{prasad-acta}
P.~Prasad.
\newblock Kinematics of a multi-dimensional shock of arbitrary strength in an
  ideal gas.
\newblock {\em Acta Mech.}, 45:163--176, 1982.

\bibitem{prasad-book}
P.~Prasad.
\newblock {\em Nonlinear hyperbolic waves in multi-dimensions}, volume 121 of
  {\em Chapman \& Hall/CRC Monographs and Surveys in Pure and Applied
  Mathematics}.
\newblock Chapman \& Hall/CRC, Boca Raton, FL, 2001.

\bibitem{prasad-parker}
P.~Prasad and D.~F. Parker.
\newblock A shock ray theory for propagation of a curved shock of arbitrary
  strength.
\newblock Technical report, Department of Mathematics, Indian Institute of
  Science, Bangalore, India, 1994.

\bibitem{sangeeta}
P.~Prasad and K.~Sangeeta.
\newblock Numerical simulation of converging nonlinear wavefronts.
\newblock {\em J. Fluid Mech.}, 385:1--20, 1999.

\bibitem{ramanathan}
T.~M. Ramanathan.
\newblock {\em Huygen's method of construction of weakly nonlinear wavefronts
  and shockfronts with application to hyperbolic caustics}.
\newblock PhD thesis, Indian Institute of Science, Bangalore, India, 1985.

\bibitem{ravindran-prasad}
R.~Ravindran and P.~Prasad.
\newblock A new theory of shock dynamics part 1: analytical considerations.
\newblock {\em Appl. Math. Lett.}, 3:77--81, 1990.

\bibitem{shu-tvd-rk}
C.-W. Shu.
\newblock Total-variation-diminishing time discretizations.
\newblock {\em SIAM J. Sci. Stat. Comput.}, 9:1073--1084, 1988.

\bibitem{srinivasan-prasad-1985}
R.~Srinivasan and P.~Prasad.
\newblock On the propagation of a multidimensional shock of arbitrary strength.
\newblock {\em Proc. Indian Acad. Sci. Math. Sci.}, 94:27--42, 1985.

\bibitem{sturtevant-kulkarny}
B.~Sturtevant and V.~A. Kulkarny.
\newblock The focusing of weak shock waves.
\newblock {\em J. Fluid Mech.}, 73:651--671, 1976.

\bibitem{takayama-etal}
K.~Takayama, H.~Kleine, and H.~Gr{\"o}nig.
\newblock An experimental investigation of the stability of converging
  cylindrical shock waves in air.
\newblock {\em Exp. Fluids}, 5:315--322, 1987.

\bibitem{whitham-book}
G.~B. Whitham.
\newblock {\em Linear and nonlinear waves}.
\newblock Wiley-Interscience [John Wiley \& Sons], New York, 1974.
\newblock Pure and Applied Mathematics.

\end{thebibliography}
\bibliographystyle{abbrv}

\end{document}